\documentclass[11pt]{article} 
\usepackage{pstricks}
\usepackage{amsmath}
\usepackage{amssymb}
\usepackage{theorem} 
\usepackage{exscale,relsize}
\usepackage{graphicx} 
\usepackage{xcolor}
\definecolor{labelkey}{rgb}{0,0.08,0.45}
\definecolor{refkey}{rgb}{0,0.6,0.0}
\definecolor{Brown}{rgb}{0.45,0.0,0.05}
\definecolor{dgreen}{rgb}{0.00,0.49,0.00}
\definecolor{dblue}{rgb}{0,0.08,0.75}
\RequirePackage[dvipdfm,colorlinks,hyperindex]{hyperref}
\hypersetup{linktocpage=true,citecolor=dblue,linkcolor=dgreen}
\title{\sffamily A PROXIMAL DECOMPOSITION METHOD FOR SOLVING \\
CONVEX VARIATIONAL INVERSE PROBLEMS\footnote{Contact author: 
P. L. Combettes, {\ttfamily plc@math.jussieu.fr},
phone: +33 1 4427 6319, fax: +33 1 4427 7200.}}
\author{Patrick L. Combettes$^1$ and Jean-Christophe Pesquet$^2$\\[5mm]
$\!^1$UPMC Universit\'e Paris 06\\
Laboratoire Jacques-Louis Lions -- UMR 7598\\
75005 Paris, France\\
{\ttfamily plc@math.jussieu.fr}\\[4mm]
$\!^2$Universit\'e Paris-Est\\
Institut Gaspard Monge and UMR CNRS 8049\\
77454 Marne la Vall\'ee Cedex 2, France\\
{\ttfamily jean-christophe.pesquet@univ-paris-est.fr}}
\date{}%\ttfamily May 6, 2009 -- version 1.12}
\newcommand{\Frac}[2]{\displaystyle{\frac{#1}{#2}}} 
\newcommand{\pscal}[2]{\langle\langle{#1}\mid{#2}\rangle\rangle} 
\newcommand{\scal}[2]{{\left\langle{{#1}\mid{#2}}\right\rangle}}
\newcommand{\menge}[2]{\big\{{#1}~\big |~{#2}\big\}} 
\newcommand{\HHH}{\ensuremath{\boldsymbol{\mathcal H}}}
\newcommand{\HH}{\ensuremath{{\mathcal H}}}
\newcommand{\GG}{\ensuremath{{\mathcal G}}}
\newcommand{\emp}{\ensuremath{{\varnothing}}}

\newcommand{\Id}{\ensuremath{\operatorname{Id}}\,}
\newcommand{\RR}{\ensuremath{\mathbb{R}}}
\newcommand{\RP}{\ensuremath{\left[0,+\infty\right[}}

\newcommand{\RPP}{\ensuremath{\left]0,+\infty\right[}}

\newcommand{\RX}{\ensuremath{\left]-\infty,+\infty\right]}}

\newcommand{\NN}{\ensuremath{\mathbb N}}
\newcommand{\KK}{\ensuremath{\mathbb K}}
\newcommand{\LL}{\ensuremath{\mathbb L}}
\newcommand{\ZZ}{\ensuremath{\mathbb Z}}
\newcommand{\weakly}{\ensuremath{\:\rightharpoonup\:}}

\newcommand{\pinf}{\ensuremath{{+\infty}}}

\newcommand{\dom}{\ensuremath{\operatorname{dom}}}
\newcommand{\prox}{\ensuremath{\operatorname{prox}}}
\newcommand{\spc}{\ensuremath{\overline{\operatorname{span}}\,}}
\newcommand{\spa}{\ensuremath{\operatorname{span}\,}}
\newcommand{\card}{\ensuremath{\operatorname{card}}}

\newcommand{\inte}{\ensuremath{\operatorname{int}}}

\newcommand{\cone}{\ensuremath{\operatorname{cone}}}
\newcommand{\sri}{\ensuremath{\operatorname{sri}}}
\newcommand{\reli}{\ensuremath{\operatorname{ri}}}

\newtheorem{theorem}{Theorem}[section]
\newtheorem{lemma}[theorem]{Lemma}

\newtheorem{proposition}[theorem]{Proposition}

\theoremstyle{plain}{\theorembodyfont{\rmfamily}%
\newtheorem{example}[theorem]{Example}}
\theoremstyle{plain}{\theorembodyfont{\rmfamily}%
\newtheorem{remark}[theorem]{Remark}}
\theoremstyle{plain}{\theorembodyfont{\rmfamily}%
\newtheorem{algorithm}[theorem]{Algorithm}}
\theoremstyle{plain}{\theorembodyfont{\rmfamily}%
}
\theoremstyle{plain}{\theorembodyfont{\rmfamily}%
}
\theoremstyle{plain}{\theorembodyfont{\rmfamily}%
}

\numberwithin{equation}{section}
\begin{document}
\maketitle

\begin{abstract}
A broad range of inverse problems can be abstracted into the problem of
minimizing the sum of several convex functions in a Hilbert space.
We propose a proximal decomposition algorithm for solving this problem 
with an arbitrary number of nonsmooth functions and establish its weak
convergence. The algorithm fully decomposes the problem in that it involves
each function individually via its own proximity operator. A significant 
improvement over the methods currently in use in the area of inverse problems 
is that it is not limited to two nonsmooth functions. Numerical applications 
to signal and image processing problems are demonstrated.
\end{abstract}

\newpage
\section{Introduction}
Throughout this paper, $\HH$ is a real Hilbert space with scalar product 
$\scal{\cdot}{\cdot}$, norm $\|\cdot\|$, and distance $d$. Moreover,
$(f_i)_{1\leq i\leq m}$ are proper lower semicontinuous convex functions 
from $\HH$ to $\RX$. We consider inverse problems that can be formulated as 
decomposed optimization problems of the form
\begin{equation}
\label{e:prob1}
\underset{x\in\HH}{\mathrm{minimize}}\;\;\sum_{i=1}^mf_i(x).
\end{equation}
In this flexible variational formulation, each potential function $f_i$ 
may represent a prior constraint on the ideal solution $\overline{x}$ 
or on the data acquisition model. The purpose of this paper is to propose 
a decomposition method that, under rather general conditions, will provide 
solutions to \eqref{e:prob1}.

To place our investigation in perspective, let us review some important 
special cases of \eqref{e:prob1} for which globally convergent numerical
methods are available. These examples encompass a variety of inverse 
problems in areas such as 
signal denoising \cite{Jsts07,Rudi92},
signal deconvolution \cite{Sign94},
Bayesian image recovery \cite{Chau07},
intensity-modulated radiation therapy \cite{Byrn05,Cens06},
image restoration \cite{Bect04,Biou07,Cham97},
linear inverse problems with sparsity constraints
\cite{Siop07,Daub07,Forn07,Trop06},
signal reconstruction from Fourier phase information \cite{Levi83},
and tomographic reconstruction \cite{Artz79,Byrn05,Star87}.

\begin{itemize}
\item[(a)]
If the functions $(f_i)_{1\leq i\leq m}$ are the indicator 
functions (see \eqref{e:iota}) of closed convex sets 
$(C_i)_{1\leq i\leq m}$ in $\HH$, \eqref{e:prob1} reduces to 
the convex feasibility problem \cite{Byrn05,Cens06,Aiep96,Star87,Youl82}
\begin{equation}
\label{e:prob11}
\text{find}\;x\in\bigcap_{i=1}^mC_i,
\end{equation}
which can be solved by projection techniques, e.g., 
\cite{Baus96,Cens84,Imag97,Kiw97b}.
\item[(b)]
The constraint sets in (a) are based on information or
measurements that can be inaccurate. As a result, the feasibility set 
$\bigcap_{i=1}^mC_i$ may turn out to be empty. An approximate solution
can be obtained by setting, for every $i\in\{1,\ldots,m\}$,
$f_i=\omega_id_{C_i}^2$, where $d_{C_i}$ is the distance function to $C_i$
(see \eqref{e:d_C}) and where $\omega_i\in\left]0,1\right]$. Thus, 
\eqref{e:prob1} becomes 
\begin{equation}
\label{e:prob12}
\underset{x\in\HH}{\mathrm{minimize}}\;\;\sum_{i=1}^m\omega_id^2_{C_i}(x).
\end{equation}
This approach is proposed in \cite{Sign94}, 
where it is solved by a parallel projection method. Finite-dimensional 
variants based on Bregman distances are investigated in \cite{Byrn01}.
\item[(c)]
If the functions $(f_i)_{1\leq i\leq m-1}$ are the indicator functions of 
closed convex sets $(C_i)_{1\leq i\leq m-1}$ in $\HH$ and 
$f_m\colon x\mapsto\|x-r\|^2$ for some $r\in\HH$, then \eqref{e:prob1} 
reduces to the best approximation problem \cite{Artz79,Sign03}
\begin{equation}
\label{e:prob14}
\underset{x\,\in\overset{m-1}{\underset{i=1}{\bigcap}}C_i}
{\mathrm{minimize}}\;\;\|x-r\|^2.
\end{equation}
Several algorithms are available to solve this problem 
\cite{Boyl86,Sign03,Gaff89,Hirs06,Yama98}. There are also methods that
are applicable in the presence of a more general strictly convex
potential $f_m$; see \cite{Sico00} and the references therein.
\item[(d)]
In \cite{Smms05}, the special instance of \eqref{e:prob1} in which $m=2$ and 
$f_2$ is Lipschitz-differentiable on $\HH$ is shown to cover a variety of
seemingly unrelated inverse problem formulations such as 
Fourier regularization problems, constrained least-squares problems, 
split feasibility problems, multiresolution sparse regularization problems, 
geometry/texture image decomposition problems, hard-constrained inconsistent 
feasibility problems, as well as certain maximum \emph{a posteriori} problems
(see also \cite{Biou07,Bred08,Lore08,Chau07,Siop07,Daub07,Forn07} for 
further developments within this framework). The forward-backward splitting 
algorithm proposed in \cite{Smms05} is governed by the updating rule
\begin{equation}
\label{e:main2005}
x_{n+1}=x_n+\lambda_n\Big(\prox_{\gamma_n f_1}
\big(x_n-\gamma_n(\nabla f_2(x_n)+b_n)\big)+a_n-x_n\Big),
\end{equation}
where $\lambda_n\in\left]0,1\right]$ and $\gamma_n\in\RPP$, where
\begin{equation}
\prox_{\gamma_nf_1}\colon x\mapsto\underset{y\in\HH}{\mathrm{argmin}}
\;\;\gamma_nf_1(y)+\frac{1}{2}\|x-y\|^2 
\end{equation}
is the proximity operator of $\gamma_nf_1$, and where the vectors 
$a_n$ and $b_n$ model tolerances in the implementation of
$\prox_{\gamma_n f_1}$ and $\nabla f_2$, respectively.
Naturally, this 2-function framework can be extended to \eqref{e:prob1} 
under the severe restriction that the functions $(f_i)_{2\leq i\leq m}$ 
be Lipschitz-differentiable. Indeed, in this case, 
$\widetilde{f_2}=\sum_{i=2}^mf_i$ also enjoys this property and 
it can be used in lieu of $f_2$ in \eqref{e:main2005}.
\item[(e)]
The problem considered in \cite{Jsts07} corresponds to $m=2$ in
\eqref{e:prob1}. In other words, the smoothness assumption on $f_2$ in 
(d) is relaxed. The algorithm adopted in \cite{Jsts07} is based
on the Douglas-Rachford splitting method \cite{Opti04,Lion79} and 
operates via the updating rule
\begin{equation}
\label{e:6ans}
\left\lfloor
\begin{array}{ll}
y_{n+\frac12}&\!\!\!\!=\prox_{\gamma f_2}y_n+a_n\\
y_{n+1}&\!\!\!\!=y_n+
\lambda_n\Big(\prox_{\gamma f_1}\big(2y_{n+\frac12}-y_n\big)
+b_n-y_{n+\frac12}\Big), 
\end{array}
\right.\\[2mm]
\end{equation}
where $\lambda_n\in\left]0,2\right[$ and $\gamma\in\RPP$, 
and where the vectors $a_n$ and $b_n$ model tolerances in 
the implementation of the proximity operators. Under suitable assumptions,
the sequence $(y_n)_{n\in\NN}$ converges weakly to a point $y\in\HH$ and 
$\prox_{\gamma f_2}y\in\operatorname{Argmin}f_1+f_2$. In this approach,
the smoothness assumption made on $f_2$ in (d) is replaced by the 
practical assumption that $\prox_{\gamma f_2}$ be implementable (to within
some error).
\end{itemize}

Some important scenarios are not covered by the above settings, namely
the formulations of type \eqref{e:prob1} that feature three or more
potentials, at least two of which are nonsmooth. In this paper, we 
investigate a reformulation of
\eqref{e:6ans} in a product space that allows us to capture instances 
of \eqref{e:prob1} in which none of the functions need be differentiable.
The resulting algorithm proceeds by decomposition in that each function is
involved individually via its own proximity operator. Since proximity
operators can be implemented for a wide variety of potentials, the proposed
framework is applicable to a broad array of problems.

In section~\ref{sec:2}, we set our notation and provide some 
background on convex analysis and proximity 
operators. We also obtain closed-form formulas for new examples of proximity 
operators that will be used subsequently. In section~\ref{sec:3}, we 
introduce our algorithm and prove its weak convergence.
Applications to signal and image processing 
problems are detailed in section~\ref{sec:4}, where numerical results are 
also provided. These results show that complex nonsmooth variational 
inverse problems, that were beyond the reach of the methods reviewed 
above, can be decomposed and solved efficiently within the proposed 
framework. Section~\ref{sec:5} concludes the paper with some remarks. 

\section{Notation and background}
\label{sec:2}

\subsection{Convex analysis}
We provide here some basic elements; for proofs and complements see 
\cite{Zali02} and, for the finite dimensional setting, \cite{Rock70}.

Let $C$ be a nonempty convex subset of $\HH$. The indicator function of $C$ is
\begin{equation}
\label{e:iota}
\iota_C\colon x\mapsto
\begin{cases}
0,&\text{if}\;\;x\in C;\\
\pinf,&\text{if}\;\;x\notin C,
\end{cases}
\end{equation}
its distance function is 
\begin{equation}
\label{e:d_C}
d_C\colon\HH\to\RP\colon x\mapsto\inf_{y\in C}\|x-y\|,
\end{equation}
its support function is
\begin{equation}
\label{e:support}
\sigma_C\colon\HH\to\RX\colon{u}\mapsto\sup_{x\in C}\scal{x}{u},
\end{equation}
and its conical hull is  
\begin{equation}
\label{e:cone}
\cone C=\bigcup_{\lambda>0}\menge{\lambda x}{x\in C}.
\end{equation}
Moreover, $\spa C$ denotes the span of $C$ and $\spc C$ the closure
of $\spa C$. The strong relative interior of $C$ is 
\begin{equation}
\label{e:sri}
\sri C=\menge{x\in C}{\cone (C-x)=\spc (C-x)}
\end{equation}
and its relative interior is 
\begin{equation}
\label{e:reli}
\reli C=\menge{x\in C}{\cone (C-x)=\spa (C-x)}.
\end{equation}
We have 
\begin{equation}
\label{e:palawan-mai08-1}
\inte C\subset\sri C\subset\reli C\subset C. 
\end{equation}

\begin{lemma}{\rm\cite[Section~6]{Rock70}}
\label{l:palawan-mai08}
Suppose that $\HH$ is finite-dimensional, and let $C$ and $D$ be 
convex subsets of $\HH$. Then the following hold.
\begin{enumerate}
\item
\label{l:palawan-mai08ii}
Suppose that $C\neq\emp$. Then $\sri C=\reli C\neq\emp$.
\item
\label{l:palawan-mai08iii}
$\reli(C-D)=\reli C-\reli D$.
\item
\label{l:palawan-mai08i}
Suppose that $D$ is an affine subspace and that $(\reli C)\cap D\neq\emp$.
Then $\reli(C\cap D)=(\reli C)\cap D$.
\end{enumerate}
\end{lemma}

Now let $C$ be a nonempty closed and convex subset of $\HH$. The 
projection of a point $x$ in $\HH$ onto $C$ is the unique point $P_Cx$ 
in $C$ such that $\|x-P_Cx\|=d_C(x)$. We have
\begin{equation}
\label{e:kolmogorov}
(\forall x\in\HH)(\forall p\in\HH)\quad
p=P_Cx\;\Leftrightarrow\;\big[\,p\in C\quad\text{and}\quad
(\forall y\in C)\quad\scal{y-p}{x-p}\leq 0\,\big].
\end{equation}
Moreover, $d_C$ is Fr\'echet differentiable on 
$\HH\smallsetminus C$ and
\begin{equation}
\label{e:Frechet}
(\forall x\in\HH\smallsetminus C)\quad\nabla d_C(x)=\frac{x-P_Cx}{d_C(x)}.
\end{equation}

The domain of a function $f\colon\HH\to\RX$ is 
$\dom f=\menge{x\in\HH}{f(x)<\pinf}$ and its set of global minimizers is
denoted by $\operatorname{Argmin}f$; if $f$ possesses a unique global 
minimizer, it is denoted by $\operatorname{argmin}_{y\in\HH}f(y)$.
The class of lower semicontinuous convex functions from $\HH$ to $\RX$ 
which are proper (i.e., with nonempty domain) is denoted by 
$\Gamma_0(\HH)$. Now let $f\in\Gamma_0(\HH)$. The conjugate of 
$f$ is the function $f^*\in\Gamma_0(\HH)$ defined by
$f^*\colon\HH\to\RX\colon u\mapsto\sup_{x\in\HH}\scal{x}{u}-f(x)$, and 
the subdifferential of $f$ is the set-valued operator
\begin{equation}
\label{e:subdiff}
\partial f\colon\HH\to 2^\HH\colon x\mapsto 
\menge{u\in\HH}{(\forall y\in\dom f)\;\:\scal{y-x}{u}+f(x)\leq f(y)}.
\end{equation}
We have 
\begin{equation}
\label{e:Fermat}
(\forall {x}\in\HH)\quad{x}\in\operatorname{Argmin}f\;
\Leftrightarrow\; 0\in\partial f(x)
\end{equation}
and
\begin{equation}
\label{e:Fenchel-Young}
(\forall x\in\HH)(\forall u\in\HH)\quad
\begin{cases}
f(x)+f^*(u)\geq\scal{x}{u}\\
f(x)+f^*(u)=\scal{x}{u}\;\Leftrightarrow\;u\in\partial f(x).
\end{cases}
\end{equation}
Moreover, if $f$ is G\^ateaux-differentiable at $x\in\HH$, then 
$\partial f(x)=\{\nabla f(x)\}$.

\begin{lemma}
\label{l:mai1968}
Let $C$ be a nonempty closed convex subset of $\HH$, let
$\phi\colon\RR\to\RR$ be an even convex function, and set 
$f=\phi\circ d_C$. Then $f\in\Gamma_0(\HH)$ and
$f^*=\sigma_C+\phi^*\circ\|\cdot\|$.
\end{lemma}
\begin{proof}
Since $\phi\colon\RR\to\RR$ is convex and even, it is continuous and 
increasing on $\RP$. On the other hand, since $C$ is convex, $d_C$ is convex. 
Hence, $\phi\,\circ d_C$ is a finite continuous convex function, which
shows that $f\in\Gamma_0(\HH)$. Moreover, 
$\phi\,\circ d_C=\phi(\inf_{y\in C}\|\cdot-y\|)
=\inf_{y\in C}\phi\circ\|\cdot-y\|$. Therefore, 
\begin{align}
(\forall u\in\HH)\quad f^*(u)
&=\sup_{x\in\HH}\scal{x}{u}-\inf_{y\in C}\phi(\|x-y\|)\nonumber\\
&=\sup_{y\in C}\scal{y}{u}+\sup_{x\in\HH}\scal{x-y}{u}
-(\phi\circ\|\cdot\|)(x-y)\nonumber\\
&=\sup_{y\in C}\scal{y}{u}+(\phi\circ\|\cdot\|)^*(u)\nonumber\\
&=\sigma_{C}(u)+(\phi\circ\|\cdot\|)^*(u).
\end{align}
Since $(\phi\,\circ\|\cdot\|)^*=\phi^*\circ\|\cdot\|$ 
\cite[Proposition~I.4.2]{Ekel99}, the proof is complete. 
\end{proof}

\subsection{Proximity operators}
For detailed accounts of the theory of proximity operators, see
\cite[Section~2]{Smms05} and \cite{More65}.

The proximity operator of a function $f\in\Gamma_0(\HH)$ is the 
operator $\prox_f\colon\HH\to\HH$ which maps every $x\in\HH$ to the 
unique minimizer of the function $f+\|x-\cdot\|^2/2$, i.e.,
\begin{equation}
\label{e:prox}
(\forall x\in\HH)\quad\prox_{f}x=\underset{y\in\HH}{\mathrm{argmin}}\;\;
f(y)+\frac12\|x-y\|^2.
\end{equation}
We have
\begin{equation}
\label{e:jjm}
(\forall x\in\HH)(\forall p\in\HH)\quad
p=\prox_{f}x\;\Leftrightarrow\;x-p\in\partial f(p).
\end{equation}
In other words, $\prox_f=(\Id+\partial f)^{-1}$.

\begin{lemma}
\label{l:0}
Let $f\in\Gamma_0(\HH)$. Then the following hold.
\begin{enumerate}
\item
\label{l:02}
$(\forall {x}\in\HH)(\forall {y}\in\HH)$ 
$\|\prox_f {x}-\prox_f {y}\|^2\leq
\scal{x-y}{\prox_f {x}-\prox_f {y}}$.
\item
\label{l:03}
$(\forall {x}\in\HH)(\forall\gamma\in\RPP)$
${x}=\prox_{\gamma{f}}{x}+
\gamma\prox_{f^*/\gamma}({x}/\gamma)$.
\end{enumerate}
\end{lemma}

\begin{lemma}{\rm\cite[Proposition~11]{Jsts07}} 
\label{l:1938}
Let $\GG$ be a real Hilbert space, let $f\in\Gamma_0(\GG)$, and let 
$L\colon\HH\to\GG$ be a bounded linear operator such that
$L\circ L^*=\kappa\Id$, for some $\kappa\in\RPP$. Then 
$f\circ L\in\Gamma_0(\HH)$ and
\begin{equation}
\label{e:pfL}
\prox_{f\circ L}=\Id+\frac{1}{\kappa}L^*\circ(\prox_{\kappa f}-\Id)\circ L.
\end{equation}
\end{lemma}

\subsection{Examples of proximity operators}

Closed-form formulas for various proximity operators are provided in 
\cite{Chau07,Siop07,Jsts07,Smms05,More65}. The following
examples will be of immediate use subsequently.

\begin{proposition}{\rm\cite[Proposition~2.10 and Remark~3.2(ii)]{Chau07}}
\label{p:decomp}
Set 
\begin{equation}
f\colon\HH\to\RX\colon x\mapsto\sum_{k\in\KK}\phi_k(\scal{x}{e_k}), 
\end{equation}
where:
\begin{enumerate}
\item
\label{p:quezoni}
$\emp\neq\KK\subset\NN$;
\item
\label{p:quezonii}
$(e_k)_{k\in\KK}$ is an orthonormal basis of $\HH$; 
\item
\label{p:quezoniii}
$(\phi_k)_{k\in\KK}$ are functions in $\Gamma_0(\RR)$;
\item
\label{p:quezoniv}
Either $\KK$ is finite, or there exists a subset $\LL$ of $\KK$ 
such that:
\begin{enumerate}
\item
\label{p:quezoniva}
$\KK\smallsetminus \LL$ is finite;
\item
\label{p:quezonivb}
$(\forall k\in\LL)$ $\phi_k\geq \phi_k(0)=0$.
\end{enumerate}
\end{enumerate}
Then $f\in\Gamma_0(\HH)$ and
\begin{equation}
(\forall x\in\HH)\quad\prox_f{x}=
\sum_{k\in\KK}\big(\prox_{\phi_k}\scal{x}{e_k}\big)e_k.
\end{equation}
\end{proposition}

We shall also require the following results, which appear to be new.

\begin{proposition}
\label{p:30avril2008}
Let $(\GG,\|\cdot\|)$ be a real Hilbert space, let $L\colon\HH\to\GG$ be 
linear and bounded, let $z\in\GG$, let $\gamma\in\RPP$, and set 
$f=\gamma\|L\cdot-z\|^2/2$. Then $f\in\Gamma_0(\HH)$ and 
\begin{equation}
\label{e:30avril2008}
(\forall x\in\HH)\quad \prox_{f}x=(\Id+\gamma L^*L)^{-1}(x+\gamma L^*z).
\end{equation}
\end{proposition}
\begin{proof}
It is clear that $f$ is a finite continuous convex function.
Now, take $x$ and $p$ in $\HH$. Then \eqref{e:jjm} yields 
$p=\prox_fx$ $\Leftrightarrow$ 
$x-p=\nabla\big(\gamma\|L\cdot-z\|^2/2\big)(p)$
$\Leftrightarrow$ $x-p=\gamma L^*(Lp-z)$ 
$\Leftrightarrow$ $p=(\Id+\gamma L^*L)^{-1}(x+\gamma L^*z)$.
\end{proof}

\begin{proposition}
\label{p:4}
Let $C$ be a nonempty closed convex subset of $\HH$, let
$\phi\colon\RR\to\RR$ be an even convex function which is 
differentiable on $\RR\smallsetminus\{0\}$, and set 
$f=\phi\circ d_C$. Then
\begin{equation}
\label{e:30octobre2007}
(\forall x\in\HH)\quad\prox_{f}x=
\begin{cases}
x+\displaystyle{\frac{\prox_{\phi^*}d_C(x)}{d_C(x)}}
(P_Cx-x),&\text{if}\;\;d_C(x)>\max\partial\phi(0);\\
P_Cx,&\text{if}\;\;d_C(x)\leq\max\partial\phi(0).
\end{cases}
\end{equation}
\end{proposition}
\begin{proof}
As seen in Lemma~\ref{l:mai1968}, $f\in\Gamma_0(\HH)$.
Now let $x\in\HH$ and set $p=\prox_f x$. Since $\phi$ is a finite even convex
function, $\partial\phi(0)=[-\beta,\beta]$ for some $\beta\in\RP$
\cite[Theorem~23.4]{Rock70}. We consider two alternatives.
\begin{itemize}
\item[(a)] 
$p\in C$: Let $y\in C$. Then $f(y)=\phi(d_C(y))=\phi(0)$
and, in particular, $f(p)=\phi(0)$. Hence, it follows from \eqref{e:jjm} and 
\eqref{e:subdiff} that 
\begin{equation}
\label{e:tashkent3mai2008}
\scal{y-p}{x-p}+\phi(0)=\scal{y-p}{x-p}+f(p)\leq f(y)=\phi(0).
\end{equation}
Consequently, $\scal{y-p}{x-p}\leq 0$ and, in view of \eqref{e:kolmogorov}, 
we get $p=P_Cx$. Thus,
\begin{equation}
\label{e:enfin2}
p\in C\quad\Leftrightarrow\quad p=P_Cx.
\end{equation}
Now, let $u\in\partial f(p)$. Since $p\in C$, $d_C(p)=0$ and, 
by \eqref{e:support}, $\sigma_C(u)\geq\scal{p}{u}$. Hence, 
\eqref{e:Fenchel-Young} and Lemma~\ref{l:mai1968} yield
\begin{equation}
\label{e:26juiin2008}
-0\,\|u\|=0\leq\sigma_C(u)-\scal{p}{u}
=\sigma_C(u)-f(p)-f^*(u)=-\phi(0)-\phi^*(\|u\|).
\end{equation}
We therefore deduce from \eqref{e:Fenchel-Young} that 
$\|u\|\in\partial\phi(0)$. Thus, $u\in\partial f(p)$ 
$\Rightarrow$ $\|u\|\leq\beta$. Since \eqref{e:jjm} 
asserts that $x-p\in\partial f(p)$, we obtain
$\|x-p\|\leq\beta$ and hence, since $p\in C$,
$d_C(x)\leq\|x-p\|\leq\beta$. As a result,
\begin{equation}
\label{e:3}
p\in C\quad\Rightarrow\quad d_C(x)\leq\beta.
\end{equation}

\item[(b)] 
$p\notin C$: Since $C$ is closed, $d_C(p)>0$ and $\phi$ is therefore 
differentiable at $d_C(p)$. It follows from \eqref{e:jjm}, the Fr\'echet 
chain rule, and \eqref{e:Frechet} that 
\begin{equation}
\label{e:30octobre2007-2}
x-p=f'(p)=\frac{\phi'(d_C(p))}{d_C(p)}(p-P_Cp).
\end{equation}
Since $\phi'\geq 0$ on $\RPP$, upon taking the norm, we obtain  
\begin{equation}
\label{e:30octobre2007-3}
\|p-x\|=\phi'(d_C(p))
\end{equation}
and therefore
\begin{equation}
\label{e:30octobre2007-4}
p-x=\frac{\|p-x\|}{d_C(p)}(P_Cp-p).
\end{equation}
In turn, appealing to Lemma~\ref{l:0}\ref{l:02} (with $f=\iota_C$) and 
\eqref{e:kolmogorov}, we obtain
\begin{equation}
\|P_Cp-P_Cx\|^2\leq\scal{p-x}{P_Cp-P_Cx}=
\frac{\|p-x\|}{d_C(p)}\scal{P_Cp-p}{P_Cp-P_Cx}\leq 0,
\end{equation}
from which we deduce that 
\begin{equation}
\label{e:30octobre2007-15}
P_Cp=P_Cx.
\end{equation}
Hence, \eqref{e:30octobre2007-4} becomes
\begin{equation}
\label{e:30octobre2007-5}
p-x=\frac{\|p-x\|}{\|p-P_Cx\|}(P_Cx-p),
\end{equation}
which can be rewritten as
\begin{equation}
\label{e:30octobre2007-6}
p-x=\frac{\|p-x\|}{\|p-x\|+\|p-P_Cx\|}(P_Cx-x).
\end{equation}
Taking the norm yields
\begin{equation}
\label{e:30octobre2007-7}
\|p-x\|=\frac{\|p-x\|}{\|p-x\|+\|p-P_Cx\|}d_C(x),
\end{equation}
and it follows from \eqref{e:30octobre2007-15} that 
\begin{equation}
\label{e:30octobre2007-7+}
d_C(x)=\|p-x\|+\|p-P_Cx\|=\|p-x\|+d_C(p).
\end{equation}
Therefore, in the light of \eqref{e:30octobre2007-3}, we obtain
\begin{equation}
\label{e:30octobre2007-11}
d_C(x)-d_C(p)=\|p-x\|=\phi'(d_C(p))
\end{equation}
and we derive from \eqref{e:jjm} that 
\begin{equation}
\label{e:30octobre2007-17}
d_C(p)=\prox_\phi d_C(x).
\end{equation}
Thus, Lemma~\ref{l:0}\ref{l:03} yields
\begin{equation}
\label{e:30octobre2007-18}
d_C(x)-d_C(p)=d_C(x)-\prox_\phi d_C(x)=\prox_{\phi^*}d_C(x)
\end{equation}
and, in turn, \eqref{e:30octobre2007-11} results in 
\begin{equation}
\label{e:30octobre2007-27}
\|p-x\|=d_C(x)-d_C(p)=\prox_{\phi^*}d_C(x).
\end{equation}
To sum up, coming back to \eqref{e:30octobre2007-6} and invoking 
\eqref{e:30octobre2007-7+} and \eqref{e:30octobre2007-27}, we obtain
\begin{align}
\label{e:30octobre2007-8}
p\notin C\quad\Rightarrow\quad p
&=x+\frac{\|p-x\|}{\|p-x\|+\|p-P_Cx\|}(P_Cx-x)\nonumber\\
&=x+\frac{\prox_{\phi^*}d_C(x)}{d_C(x)}(P_Cx-x).
\end{align}
Furthermore, we derive from \eqref{e:30octobre2007-17} and \eqref{e:jjm}
that 
\begin{equation}
\label{e:4}
p\notin C 
\;\Rightarrow\;d_C(p)>0
\;\Rightarrow\;\prox_\phi d_C(x)\neq 0
\;\Rightarrow\;d_C(x)\notin\partial\phi(0)
\;\Rightarrow\;d_C(x)>\beta. 
\end{equation}
\end{itemize}
Upon combining \eqref{e:3} and \eqref{e:4}, we obtain
\begin{equation}
\label{e:enfin1}
p\in C\quad\Leftrightarrow\quad d_C(x)\leq\beta.
\end{equation}
Altogether, \eqref{e:30octobre2007} follows from \eqref{e:enfin2}, 
\eqref{e:30octobre2007-8}, and \eqref{e:enfin1}.
\end{proof}

The above proposition shows that a nice feature of the proximity operator
of $\phi\circ d_C$ is that it can be decomposed in terms of 
$\prox_{\phi^*}$ and $P_C$. Here is an application of this result.
\begin{proposition}
\label{p:manille2008-05-05}
Let $C$ be a nonempty closed convex subset of $\HH$, let
$\alpha\in\RPP$, let $p\in\left[1,\pinf\right[$, and set 
$f=\alpha d^p_C$. Then the following hold.
\begin{enumerate}
\item
\label{p:manille2008-05-05i}
Suppose that $p=1$. Then 
\begin{equation}
\label{e:6mai2008-1}
(\forall x\in\HH)\quad\prox_{f}x=
\begin{cases}
x+\displaystyle{\frac{\alpha}{d_C(x)}}
(P_Cx-x),&\text{if}\;\;d_C(x)>\alpha;\\
P_Cx,&\text{if}\;\;d_C(x)\leq\alpha.
\end{cases}
\end{equation}
\item
\label{p:manille2008-05-05ii}
Suppose that $p>1$. Then 
\begin{equation}
\label{e:5mai2008-1}
(\forall x\in\HH)\quad\prox_{f}x=
\begin{cases}
x+\displaystyle{\frac{\nu(x)}{d_C(x)}}
(P_Cx-x),&\text{if}\;\;x\notin C;\\
x,&\text{if}\;\;x\in C,
\end{cases}
\end{equation}
where $\nu(x)$ is the unique real number in $\RP$ that satisfies
$\nu(x)+(\nu(x)/(\alpha p))^{1/(p-1)}=d_C(x)$.
\end{enumerate}
\end{proposition}
\begin{proof}
\ref{p:manille2008-05-05i}:
Set $\phi=\alpha|\cdot|$. Then
$\max\partial\phi(0)=\max\,[-\alpha,\alpha]=\alpha$ and
$\phi^*=\iota_{[-\alpha,\alpha]}$. Therefore, 
$\prox_{\phi^*}=P_{[-\alpha,\alpha]}$ and hence
$(\forall\mu\in\left]\alpha,\pinf\right[)$ 
$\prox_{\phi^*}\mu=\alpha$. In view of \eqref{e:30octobre2007}, we obtain 
\eqref{e:6mai2008-1}.

\ref{p:manille2008-05-05ii}:
Let $x\in\HH$ and note that, since $C$ is closed, $d_C(x)>0$
$\Leftrightarrow$ $x\notin C$. Now set $\phi=\alpha|\cdot|^p$. 
Then $\max\partial\phi(0)=\max\{0\}=0$ and
$\phi^*\colon\mu\mapsto (p-1)(\alpha p)^{1/(1-p)}|\mu|^{p/(p-1)}/p$. 
Hence, it follows from \eqref{e:jjm} and 
\cite[Corollary~2.5]{Siop07} that $\prox_{\phi^*}d_C(x)$
is the unique solution $\nu(x)\in\RP$ to the equation
$d_C(x)-\nu(x)=\phi^{*\,\prime}(\nu(x))=
(\nu(x)/(\alpha p))^{1/(p-1)}$.
Appealing to \eqref{e:30octobre2007}, we obtain \eqref{e:5mai2008-1}.
\end{proof}

Let us note that explicit expressions can be obtained for several values of
$p$ in Proposition~\ref{p:manille2008-05-05}\ref{p:manille2008-05-05ii}. 
Here is an example that will be used subsequently.

\begin{example}
\label{ex:manila2008-05-07}
Let $C$ be a nonempty closed convex subset of $\HH$, let
$\alpha\in\RPP$, and set $f=\alpha d^{3/2}_C$. Then 
\begin{equation}
\label{e:7mai2008}
(\forall x\in\HH)\quad\prox_{f}x=
\begin{cases}
x+\displaystyle{\frac{9\alpha^2\big(\sqrt{1+16d_C(x)/(9\alpha^2)}-1\big)}
{8d_C(x)}}(P_Cx-x),&\text{if}\;\;x\notin C;\\
x,&\text{if}\;\;x\in C.
\end{cases}
\end{equation}
\end{example}
\begin{proof}
Set $p=3/2$ in
Proposition~\ref{p:manille2008-05-05}\ref{p:manille2008-05-05ii}.
\end{proof}

\section{Algorithm and convergence}
\label{sec:3}
The main algorithm is presented in section~\ref{sec:31}.
In section~\ref{sec:32}, we revisit the Douglas-Rachford algorithm in the
context of minimization problems (Proposition~\ref{p:2}), with special 
emphasis on its convergence in a specific case (Proposition~\ref{p:3}). 
These results are transcribed in a product space in section~\ref{sec:33} 
to prove the weak convergence of Algorithm~\ref{algo:1}.

\subsection{Algorithm}
\label{sec:31}
We propose the following proximal method to solve \eqref{e:prob1}.
In this splitting algorithm, each function $f_i$ is used separately by
means of its own proximity operator. 

\begin{algorithm}
\label{algo:1}
For every $i\in\{1,\ldots,m\}$, let $(a_{i,n})_{n\in\NN}$ be a sequence 
in $\HH$. A sequence $(x_n)_{n\in\NN}$ is generated by the 
following routine.
\begin{equation}
\label{e:main1}
\begin{array}{l}
\text{Initialization}\\
\left\lfloor
\begin{array}{l}
\gamma\in\RPP\\[1mm]
(\omega_{i})_{1\leq i\leq m}\in\left]0,1\right]^m
\;\text{satisfy}\;\sum_{i=1}^m\omega_i=1\\[1mm]
(y_{i,0})_{1\leq i\leq m}\in\HH^m\\[1mm]
x_0=\displaystyle{\sum_{i=1}^m}\,\omega_iy_{i,0}\\[3mm]
\end{array}
\right.\\[13mm]
\text{For}\;n=0,1,\ldots\\
\left\lfloor
\begin{array}{l}
\text{For}\;i=1,\ldots,m\\
\quad\left\lfloor
\begin{array}{l}
p_{i,n}=\prox_{\gamma f_i/\omega_i}y_{i,n}+a_{i,n}\\
\end{array}
\right.\\[2mm]
p_n=\displaystyle{\sum_{i=1}^m}\,\omega_ip_{i,n}\\[5mm]
\lambda_n\in\left]0,2\right[\\[2mm]
\text{For}\;i=1,\ldots,m\\
\quad\left\lfloor
\begin{array}{l}
y_{i,n+1}=y_{i,n}+\lambda_n\big(2p_n-x_n-p_{i,n}\big)
\end{array}
\right.\\[2mm]
x_{n+1}=x_n+\lambda_n(p_{n}-x_n).
\end{array}
\right.\\
\end{array}
\end{equation}
\end{algorithm}

At iteration $n$, the proximal vectors $(p_{i,n})_{1\leq i\leq m}$, as well
as the auxiliary vectors $(y_{i,n})_{1\leq i\leq m}$, can be computed 
simultaneously, hence the parallel structure of Algorithm~\ref{algo:1}.
Another feature of the algorithm is that some error $a_{i,n}$ is tolerated
in the computation of the $i$th proximity operator. 

\subsection{The Douglas-Rachford algorithm for minimization problems}
\label{sec:32}
To ease our presentation, we introduce in this section a second real 
Hilbert space $(\HHH,|||\cdot|||)$. As usual, $\weakly$ denotes weak
convergence.

The (nonlinear) Douglas-Rachford splitting method was initially developed 
for the problem of finding a zero of the sum of two maximal monotone operators
in \cite{Lion79} (see \cite{Opti04} for recent refinements). In 
the case when the maximal monotone operators are 
subdifferentials, it provides an algorithm for minimizing the sum of two
convex functions. In this section, we develop this point of view, starting
with the following result.

\begin{proposition}
\label{p:2}
Let ${\boldsymbol f}_1$ and ${\boldsymbol f}_2$ be functions in 
$\Gamma_0(\HHH)$, let $({\boldsymbol a}_n)_{n\in\NN}$ and
$({\boldsymbol b}_n)_{n\in\NN}$ be sequences in $\HHH$, and let 
$({\boldsymbol y}_n)_{n\in\NN}$ be a sequence generated by the 
following routine.
\begin{equation}
\label{e:main23}
\begin{array}{l}
\mathrm{Initialization}\\
\left\lfloor
\begin{array}{l}
\gamma\in\RPP\\
{\boldsymbol y}_0\in\HHH\\[1mm]
\end{array}
\right.\\[5mm]
\mathrm{For}\;n=0,1,\ldots\\
\left\lfloor
\begin{array}{l}
{\boldsymbol y}_{n+\frac12}=\prox_{\gamma {\boldsymbol f}_2}{\boldsymbol y}_n
+{\boldsymbol a}_n\\[2mm]
\lambda_n\in\left]0,2\right[\\
{\boldsymbol y}_{n+1}={\boldsymbol y}_n+
\lambda_n\Big(\prox_{\gamma {\boldsymbol f}_1}\big(2{\boldsymbol y}_{n+\frac12}-
{\boldsymbol y}_n\big)+{\boldsymbol b}_n-{\boldsymbol y}_{n+\frac12}\Big). 
\end{array}
\right.
\end{array}
\end{equation}
Set 
\begin{equation}
\label{e:nadal1}
{\boldsymbol G}=\operatorname{Argmin}{\boldsymbol f}_1+{\boldsymbol f}_2
\quad\text{and}\quad 
{\boldsymbol T}=2\prox_{\gamma{\boldsymbol f}_1}\circ\,
(2\prox_{\gamma{\boldsymbol f}_2}-\,{\boldsymbol\Id})-
2\prox_{\gamma {\boldsymbol f}_2}+\,{\boldsymbol\Id},
\end{equation}
and suppose that the following hold.
\begin{enumerate}
\item
\label{p:2i}
$\lim\limits_{|||{\boldsymbol x}|||\to\pinf}{\boldsymbol f}_1({\boldsymbol x})
+{\boldsymbol f}_2({\boldsymbol x})=\pinf$.
\item
\label{p:2ii}
${\boldsymbol 0}\in\sri(\dom {\boldsymbol f}_1-\dom {\boldsymbol f}_2)$.
\item
\label{p:2iii}
$\sum_{n\in\NN}\lambda_n(2-\lambda_n)=\pinf$.
\item
\label{p:2iv}
$\sum_{n\in\NN}\lambda_n(|||{\boldsymbol a}_n|||+
|||{\boldsymbol b}_n|||)<\pinf$. 
\end{enumerate}
Then ${\boldsymbol G}\neq\emp$, $({\boldsymbol y}_n)_{n\in\NN}$ converges 
weakly to a fixed point ${\boldsymbol y}$ of ${\boldsymbol T}$, and 
$\prox_{\gamma{\boldsymbol f}_2}{\boldsymbol y}\in {\boldsymbol G}$. 
\end{proposition}
\begin{proof}
It follows from \ref{p:2ii} that 
$\dom({\boldsymbol f}_1+{\boldsymbol f}_2)=
\dom{\boldsymbol f}_1\cap\dom{\boldsymbol f}_2\neq\emp$. Hence,
since ${\boldsymbol f}_1+{\boldsymbol f}_2$ is lower semicontinuous 
and convex as the sum of two such functions, we have 
${\boldsymbol f}_1+{\boldsymbol f}_2\in\Gamma_0(\HHH)$. In turn,
we derive from \ref{p:2i} and \cite[Theorem~2.5.1(ii)]{Zali02} that 
\begin{equation}
\label{e:GG}
{\boldsymbol G}\neq\emp.
\end{equation}
Next, let us set ${\boldsymbol A}_1=\partial {\boldsymbol f}_1$,
${\boldsymbol A}_2=\partial {\boldsymbol f}_2$, and 
${\boldsymbol Z}=\menge{{\boldsymbol x}\in\HHH}{{\boldsymbol 0}
\in{\boldsymbol A}_1{\boldsymbol x}+{\boldsymbol A}_2{\boldsymbol x}}$.
Then ${\boldsymbol A}_1$ and ${\boldsymbol A}_2$ are maximal monotone 
operators \cite[Theorem~3.1.11]{Zali02}. In addition, in view of
\eqref{e:jjm}, the resolvents
of $\gamma{\boldsymbol A}_1$ and $\gamma{\boldsymbol A}_2$ are respectively
\begin{equation}
\label{e:28avril2008}
J_{\gamma{\boldsymbol A}_1}=(\Id+\gamma{\boldsymbol A}_1)^{-1}
=\prox_{\gamma{\boldsymbol f}_1}\quad\text{and}\quad
J_{\gamma{\boldsymbol A}_2}=(\Id+\gamma{\boldsymbol A}_2)^{-1}
=\prox_{\gamma{\boldsymbol f}_2}.
\end{equation}
Thus, the iteration in \eqref{e:main23} can be rewritten as
\begin{equation}
\label{e:7ans}
\left\lfloor
\begin{array}{l}
{\boldsymbol y}_{n+\frac12}=J_{\gamma {\boldsymbol A}_2}{\boldsymbol y}_n
+{\boldsymbol a}_n\\[2mm]
\lambda_n\in\left]0,2\right[\\
{\boldsymbol y}_{n+1}={\boldsymbol y}_n+\lambda_n\Big(
J_{\gamma {\boldsymbol A}_1}\big(2{\boldsymbol y}_{n+\frac12}
-{\boldsymbol {\boldsymbol y}}_n\big)
+{\boldsymbol b}_n-{\boldsymbol y}_{n+\frac12}\Big). 
\end{array}
\right.
\end{equation}
Moreover, it follows from \eqref{e:Fermat}, \ref{p:2ii}, and 
\cite[Theorem~2.8.3]{Zali02} that 
\begin{equation}
\label{e:c5}
{\boldsymbol G}=\menge{{\boldsymbol x}\in\HHH}
{{\boldsymbol 0}\in\partial({\boldsymbol f}_1+
{\boldsymbol f}_2)({\boldsymbol x})}
=\menge{{\boldsymbol x}\in\HHH}{{\boldsymbol 0}\in\partial{\boldsymbol f}_1
({\boldsymbol x})+\partial{\boldsymbol f}_2({\boldsymbol x})}
={\boldsymbol Z}.
\end{equation}
Thus, \eqref{e:GG} yields ${\boldsymbol Z}\neq\emp$ and 
it follows from \ref{p:2iii}, \ref{p:2iv}, and the results of 
\cite[Section~5]{Opti04} that $({\boldsymbol y}_n)_{n\in\NN}$ converges 
weakly to a fixed point ${\boldsymbol y}$ of the operator
$2J_{\gamma{\boldsymbol A}_1}\circ (2J_{\gamma{\boldsymbol A}_2}-
\,{\boldsymbol\Id})-2J_{\gamma {\boldsymbol A}_2}+\,{\boldsymbol\Id}$,
and that $J_{\gamma {\boldsymbol A}_2}{\boldsymbol y}\in {\boldsymbol Z}$.
In view of \eqref{e:nadal1}, \eqref{e:28avril2008}, and \eqref{e:c5}, 
the proof is complete.
\end{proof}

It is important to stress that algorithm~\eqref{e:main23} provides 
a minimizer indirectly: the sequence $({\boldsymbol y}_n)_{n\in\NN}$ 
is first constructed, and then a minimizer of 
${\boldsymbol f}_1+{\boldsymbol f}_2$ is obtained as the image of 
the weak limit ${\boldsymbol y}$ of
$({\boldsymbol y}_n)_{n\in\NN}$ under 
$\prox_{\gamma {\boldsymbol f}_2}$. In general, nothing is known
about the weak convergence of the sequences
$(\prox_{\gamma {\boldsymbol f}_1}{\boldsymbol y}_n)_{n\in\NN}$ and
$(\prox_{\gamma {\boldsymbol f}_2}{\boldsymbol y}_n)_{n\in\NN}$.
The following result describes a remarkable situation in which 
$(\prox_{\gamma {\boldsymbol f}_1}{\boldsymbol y}_n)_{n\in\NN}$ 
does converges weakly and its weak limit turns out to be a minimizer of 
${\boldsymbol f}_1+{\boldsymbol f}_2$.

\begin{proposition}
\label{p:3}
Let ${\boldsymbol D}$ be a closed vector subspace of $\HHH$, let 
${\boldsymbol f}\in\Gamma_0(\HHH)$, let
$({\boldsymbol a}_n)_{n\in\NN}$ be a sequence in $\HHH$, and let 
$({\boldsymbol x}_n)_{n\in\NN}$ be a sequence generated by the 
following routine.
\begin{equation}
\label{e:main24}
\begin{array}{l}
\mathrm{Initialization}\\
\left\lfloor
\begin{array}{l}
\gamma\in\RPP\\
{\boldsymbol y}_0\in\HHH\\
{\boldsymbol x}_0=P_{{\boldsymbol D}}\,{\boldsymbol y}_0\\[1mm]
\end{array}
\right.\\[5mm]
\mathrm{For}\;n=0,1,\ldots\\
\left\lfloor
\begin{array}{l}
{\boldsymbol y}_{n+\frac12}=\prox_{\gamma {\boldsymbol f}}{\boldsymbol y}_n
+{\boldsymbol a}_n\\
{\boldsymbol p}_{n}=P_{{\boldsymbol D}}\,{\boldsymbol y}_{n+\frac12}\\[2mm]
\lambda_n\in\left]0,2\right[\\
{\boldsymbol y}_{n+1}={\boldsymbol y}_n+\lambda_n\big(2{\boldsymbol p}_{n}-
{\boldsymbol x}_n-{\boldsymbol y}_{n+\frac12}\big)\\
{\boldsymbol x}_{n+1}={\boldsymbol x}_n+\lambda_n({\boldsymbol p}_{n}-
{\boldsymbol x}_n). 
\end{array}
\right.
\end{array}
\end{equation}
Let ${\boldsymbol G}$ be the set of minimizers of ${\boldsymbol f}$ over
${\boldsymbol D}$ and suppose that the following hold.
\begin{enumerate}
\item
\label{p:3i}
$\lim\limits_{{\boldsymbol x}\in {\boldsymbol D},\:
|||{\boldsymbol x}|||\to\pinf}\;{\boldsymbol f}({\boldsymbol x})
=\pinf$.
\item
\label{p:3ii}
${\boldsymbol 0}\in\sri({\boldsymbol D}-\dom {\boldsymbol f})$.
\item
\label{p:3iii}
$\sum_{n\in\NN}\lambda_n(2-\lambda_n)=\pinf$.
\item
\label{p:3iv}
$\sum_{n\in\NN}\lambda_n|||{\boldsymbol a}_n|||<\pinf$. 
\end{enumerate}
Then ${\boldsymbol G}\neq\emp$ and $({\boldsymbol x}_n)_{n\in\NN}$ converges 
weakly to a point in ${\boldsymbol G}$. 
\end{proposition}
\begin{proof}
Set ${\boldsymbol f}_1=\iota_{\boldsymbol D}$,
${\boldsymbol f}_2={\boldsymbol f}$,
and $(\forall n\in\NN)$ ${\boldsymbol b}_n={\boldsymbol 0}$. 
Then \eqref{e:iota} and 
\eqref{e:prox} yield $\prox_{\gamma{\boldsymbol f}_1}=P_{\boldsymbol D}$ 
and, since ${\boldsymbol D}$ is a closed vector subspace, $P_{\boldsymbol D}$ 
is a linear operator. Hence, proceeding by induction, we can rewrite 
the update equation for ${\boldsymbol x}_n$ in \eqref{e:main24} as
\begin{align}
{\boldsymbol x}_{n+1}
&={\boldsymbol x}_n+\lambda_n({\boldsymbol p}_{n}-
{\boldsymbol x}_n)\nonumber\\ 
&=P_{{\boldsymbol D}}\,{\boldsymbol y}_n+\lambda_n\big(2P_{{\boldsymbol D}}\,
{\boldsymbol p}_{n}-P_{{\boldsymbol D}}\,{\boldsymbol x}_n
-P_{{\boldsymbol D}}\,{\boldsymbol y}_{n+\frac12}\big)\nonumber\\
&=P_{{\boldsymbol D}}\,\Big({\boldsymbol y}_n+\lambda_n\big(2
{\boldsymbol p}_{n}-{\boldsymbol x}_n-
{\boldsymbol y}_{n+\frac12}\big)\Big)\nonumber\\
&=P_{\boldsymbol D}\,{\boldsymbol y}_{n+1}.
\end{align}
As a result, \eqref{e:main24} is equivalent to 
\begin{equation}
\label{e:main24-5}
\begin{array}{l}
\mathrm{Initialization}\\
\left\lfloor
\begin{array}{l}
\gamma\in\RPP\\
{\boldsymbol y}_0\in\HHH\\[1mm]
\end{array}
\right.\\[5mm]
\mathrm{For}\;n=0,1,\ldots\\
\left\lfloor
\begin{array}{l}
{\boldsymbol x}_n=P_{{\boldsymbol D}}\,{\boldsymbol y}_n\\
{\boldsymbol y}_{n+\frac12}=\prox_{\gamma {\boldsymbol f}}{\boldsymbol y}_n
+{\boldsymbol a}_n\\
{\boldsymbol p}_{n}=P_{{\boldsymbol D}}\,{\boldsymbol y}_{n+\frac12}\\[2mm]
\lambda_n\in\left]0,2\right[\\
{\boldsymbol y}_{n+1}={\boldsymbol y}_n+\lambda_n\big(2{\boldsymbol p}_{n}-
{\boldsymbol x}_n-{\boldsymbol y}_{n+\frac12}\big). 
\end{array}
\right.
\end{array}
\end{equation}
Thus, since
\begin{equation}
(\forall{\boldsymbol x}\in\HHH)(\forall{\boldsymbol y}\in\HHH)\quad
P_{\boldsymbol D}(2{\boldsymbol y}-{\boldsymbol x})=
2P_{\boldsymbol D}{\boldsymbol y}-P_{\boldsymbol D}{\boldsymbol x},
\end{equation}
\eqref{e:main24-5} appears as a special case of \eqref{e:main23} in which 
we have introduced the auxiliary variables ${\boldsymbol x}_n$ and
${\boldsymbol p}_n$. In addition, the operator ${\boldsymbol T}$ of 
\eqref{e:nadal1} becomes
\begin{equation}
\label{e:nadal2}
{\boldsymbol T}=4(P_{\boldsymbol D}\circ\prox_{\gamma{\boldsymbol f}})
-2P_{\boldsymbol D}-2\prox_{\gamma{\boldsymbol f}}+\,{\boldsymbol\Id}.
\end{equation}
Since \ref{p:3i}--\ref{p:3iv} are specializations of their
respective counterparts in Proposition~\ref{p:2}, it follows from 
Proposition~\ref{p:2} that ${\boldsymbol G}\neq\emp$ and that 
there exists a fixed point 
${\boldsymbol y}$ of ${\boldsymbol T}$ such that 
${\boldsymbol y}_n\weakly{\boldsymbol y}$ and
$\prox_{\gamma{\boldsymbol f}}{\boldsymbol y}\in {\boldsymbol G}$. 
Note that, since ${\boldsymbol G}\subset {\boldsymbol D}$, 
$\prox_{\gamma{\boldsymbol f}}{\boldsymbol y}\in{\boldsymbol D}$ and, in
turn, $P_{\boldsymbol D}(\prox_{\gamma{\boldsymbol f}}{\boldsymbol y})=
\prox_{\gamma{\boldsymbol f}}{\boldsymbol y}$.
Thus, in view of \eqref{e:nadal2}, we obtain
\begin{eqnarray}
{\boldsymbol T}{\boldsymbol y}={\boldsymbol y}
&\Leftrightarrow&
4P_{\boldsymbol D}(\prox_{\gamma{\boldsymbol f}}
{\boldsymbol y})-2P_{\boldsymbol D}{\boldsymbol y}-
2\prox_{\gamma{\boldsymbol f}}{\boldsymbol y}+{\boldsymbol y}
={\boldsymbol y}\\
&\Leftrightarrow&
2P_{\boldsymbol D}(\prox_{\gamma{\boldsymbol f}}{\boldsymbol y})-
P_{\boldsymbol D}{\boldsymbol y}=
\prox_{\gamma{\boldsymbol f}}{\boldsymbol y}\nonumber\\
&\Leftrightarrow&
\prox_{\gamma{\boldsymbol f}}{\boldsymbol y}
=P_{\boldsymbol D}{\boldsymbol y}.
\end{eqnarray}
Hence, since 
$\prox_{\gamma{\boldsymbol f}}{\boldsymbol y}\in{\boldsymbol G}$, we
also have $P_{\boldsymbol D}{\boldsymbol y}\in {\boldsymbol G}$. On 
the other hand, since $P_{\boldsymbol D}$ is linear and continuous, it is
weakly continuous and therefore ${\boldsymbol y}_n\weakly {\boldsymbol y}$ 
$\Rightarrow$ $P_{\boldsymbol D}{\boldsymbol y}_n\weakly 
{\boldsymbol P}_{\boldsymbol D}{\boldsymbol y}\in {\boldsymbol G}$. 
We conclude that ${\boldsymbol x}_n\weakly 
{\boldsymbol P}_{\boldsymbol D}{\boldsymbol y}\in {\boldsymbol G}$. 
\end{proof}

\subsection{Convergence of Algorithm~\ref{algo:1}}
\label{sec:33}

The convergence of the main algorithm can now be demonstrated.

\begin{theorem}
\label{t:3}
Let $G$ be the set of solutions to \eqref{e:prob1} and 
let $(x_n)_{n\in\NN}$ be a sequence generated by Algorithm~\ref{algo:1}
under the following assumptions.
\begin{enumerate}
\item
\label{t:3i}
$\lim\limits_{\|x\|\to\pinf}f_1(x)+\cdots+f_m(x)=\pinf$.
\item
\label{t:3ii}
$(0,\ldots,0)\in\sri\menge{(x-x_1,\ldots,x-x_m)}
{x\in\HH,\,x_1\in\dom f_1,\ldots,\, x_m\in\dom f_m}$.
\item
\label{t:3iii}
$\sum_{n\in\NN}\lambda_n(2-\lambda_n)=\pinf$.
\item
\label{t:3iv}
$(\forall i\in\{1,\ldots,m\})$
$\sum_{n\in\NN}\lambda_n\|a_{i,n}\|<\pinf$. 
\end{enumerate}
Then $G\neq\emp$ and $(x_n)_{n\in\NN}$ converges weakly to a point in $G$.
\end{theorem}
\begin{proof}
Let $\HHH$ be the real Hilbert space obtained by endowing the $m$-fold 
Cartesian product $\HH^m$ with the scalar product
\begin{equation}
\label{e:directscal}
\pscal{\cdot}{\cdot}\colon(\boldsymbol{x},\boldsymbol{y})\mapsto
\sum_{i=1}^m\omega_i\scal{x_i}{y_i}, 
\end{equation}
where $(\omega_i)_{1\leq i\leq m}$ is defined in \eqref{e:main1}, and 
where $\boldsymbol{x}=(x_i)_{1\leq i\leq m}$ and 
$\boldsymbol{y}=(y_i)_{1\leq i\leq m}$ denote generic elements in $\HHH$. 
The associated norm is denoted by $|||\cdot|||$, i.e.,
\begin{equation}
\label{e:directnorm}
|||\cdot|||\colon\boldsymbol{x}\mapsto
\sqrt{\sum_{i=1}^m\omega_i\|x_i\|^2}.
\end{equation}
Furthermore, set
\begin{equation}
\label{e:D}
\boldsymbol{D}=\menge{(x,\dots,x)\in\HHH}{x\in\HH}
\end{equation}
and
\begin{equation}
\label{e:F}
{\boldsymbol f}\colon\HHH\to\RX\colon {\boldsymbol x}
\mapsto\sum_{i=1}^mf_i(x_i).
\end{equation}
It follows from \eqref{e:directnorm} that $\boldsymbol{D}$ is 
a closed vector subspace of $\HHH$ with projector
\begin{equation}
\label{e:PD}
P_{\boldsymbol{D}}\colon\boldsymbol{x}\mapsto
\bigg(\sum_{i=1}^m\omega_ix_i,\ldots,\sum_{i=1}^m\omega_ix_i\bigg),
\end{equation}
and that the operator
\begin{equation}
\label{e:iso}
{\boldsymbol j}\colon\HH\to{\boldsymbol D}\colon x\mapsto(x,\ldots,x) 
\end{equation}
is an isomorphism. In addition, ${\boldsymbol f}\in\Gamma_0(\HHH)$ and we 
derive from \eqref{e:prox}, \eqref{e:directnorm}, and \eqref{e:F} that 
\begin{equation}
\label{e:ProxF}
\prox_{\boldsymbol{f}}\colon\boldsymbol{x}\mapsto
\big(\prox_{f_1/\omega_1}x_1,\ldots,\prox_{f_m/\omega_m}x_m\big).
\end{equation}
From the sequences $(x_n)_{n\in\NN}$, $(p_n)_{n\in\NN}$, 
$((y_{i,n})_{n\in\NN})_{1\leq i\leq m}$, 
$((p_{i,n})_{n\in\NN})_{1\leq i\leq m}$, 
and $((a_{i,n})_{n\in\NN})_{1\leq i\leq m}$ of Algorithm~\ref{algo:1} 
we define, for every $n\in\NN$,
\begin{equation}
\label{e:nadal3}
{\boldsymbol x}_n={\boldsymbol j}(x_{n}),\;\:
{\boldsymbol p}_n={\boldsymbol j}(p_{n}),\;\:
{\boldsymbol y}_n=(y_{i,n})_{1\leq i\leq m},\;\:
{\boldsymbol y}_{n+1/2}=(p_{i,n})_{1\leq i\leq m},\;\text{and}\;\:
{\boldsymbol a}_n=(a_{i,n})_{1\leq i\leq m}.
\end{equation}
It follows from \eqref{e:PD}, \eqref{e:iso}, and \eqref{e:ProxF} that 
the sequences defined in \eqref{e:nadal3} are precisely those involved 
in \eqref{e:main24}, and that the set of minimizers ${\boldsymbol G}$ in
Proposition~\ref{p:3} is precisely
\begin{equation}
\label{e:G}
{\boldsymbol G}={\boldsymbol j}(G).
\end{equation}
On the other hand,
it follows from \eqref{e:directnorm}, \eqref{e:D}, and \eqref{e:F} that 
the properties \ref{t:3i}--\ref{t:3iv} above yield their
respective counterparts in Proposition~\ref{p:3}. Thus, we deduce from
Proposition~\ref{p:3} and \eqref{e:G} that $({\boldsymbol x}_n)_{n\in\NN}$ 
converges weakly to a point ${\boldsymbol j}(x)$ for some $x\in G$.
Thus, $(x_n)_{n\in\NN}=({\boldsymbol j}^{-1}({\boldsymbol x}_n))_{n\in\NN}$ 
converges weakly to $x\in G$. 
\end{proof}

\begin{remark}\
\begin{enumerate}
\item
We have conveniently obtained Algorithm~\ref{algo:1} as a direct 
transcription of a special case (see Proposition~\ref{p:3}) of the 
Douglas-Rachford algorithm transposed in a product space. A similar 
decomposition method could be obtained by using the theory of partial 
inverses for monotone operators \cite{Spin85}.
\item
When $m=2$, Algorithm~\ref{algo:1} does not revert to the standard
Douglas-Rachford iteration \eqref{e:6ans}.
Actually, even in this case, it seems better to use the former to the
extent that, as seen in Theorem~\ref{t:3}, it produces directly a sequence 
that converges weakly to a minimizer of $f_1+f_2$.
\end{enumerate}
\end{remark}

To conclude this section, we describe some situations in which 
condition~\ref{t:3ii} in Theorem~\ref{t:3} is satisfied.
\begin{proposition}
\label{p:sabang-mai2008}
Set ${\boldsymbol C}=\menge{(x-x_1,\ldots,x-x_m)}
{x\in\HH,\,x_1\in\dom f_1,\ldots,\, x_m\in\dom f_m}$ and
suppose that any of the following holds.
\begin{enumerate}
\item
\label{p:sabang-mai2008i}
${\boldsymbol C}$ is a closed vector subspace.
\item
\label{p:sabang-mai2008ii}
$\bigcap_{i=1}^m\dom f_i\neq\emp$ and
$(\dom f_i)_{1\leq i\leq m}$ are affine subspaces of finite dimensions.
\item
\label{p:sabang-mai2008iii}
$\bigcap_{i=1}^m\dom f_i\neq\emp$ and
$(\dom f_i)_{1\leq i\leq m}$ are closed affine subspaces of 
finite codimensions.
\item
\label{p:sabang-mai2008iv}
${\boldsymbol 0}\in\inte{\boldsymbol C}$.
\item
\label{p:sabang-mai2008v}
$\dom f_1\cap\bigcap_{i=2}^m\inte\dom f_i\neq\emp$.
\item
\label{p:sabang-mai2008vi}
$\HH$ is finite-dimensional and $\bigcap_{i=1}^m\reli\dom f_i\neq\emp$.
\end{enumerate}
Then ${\boldsymbol 0}\in\sri{\boldsymbol C}$.
\end{proposition}
\begin{proof}
We use the notation of the proof of Theorem~\ref{t:3}, hence
${\boldsymbol C}={\boldsymbol D}-\dom{\boldsymbol f}$.

\ref{p:sabang-mai2008i}: 
We have ${\boldsymbol C}=\spc{\boldsymbol C}$. Since
${\boldsymbol C}\subset\cone{\boldsymbol C}\subset\spa{\boldsymbol C}
\subset\spc{\boldsymbol C}$, we therefore obtain 
$\cone {\boldsymbol C}=\spc{\boldsymbol C}$. Appealing to 
\eqref{e:sri}, we conclude that ${\boldsymbol 0}\in\sri {\boldsymbol C}$.

\ref{p:sabang-mai2008ii}$\Rightarrow$\ref{p:sabang-mai2008i}: 
The assumption implies that 
$\dom {\boldsymbol f}=\dom f_1\times\cdots\times\dom f_m$ is a 
finite-dimensional affine subspace of $\HHH$ and that
${\boldsymbol D}\cap\dom {\boldsymbol f}\neq\emp$. 
Since ${\boldsymbol D}$ is closed 
vector subspace, it follows from \cite[Lemma~9.36]{Deut01} that 
${\boldsymbol D}-\dom{\boldsymbol f}$ is a closed vector subspace.

\ref{p:sabang-mai2008iii}$\Rightarrow$\ref{p:sabang-mai2008i}: 
Here $\dom {\boldsymbol f}=\dom f_1\times\cdots\times\dom f_m$ is a closed
affine subspace of $\HHH$ of finite codimension and that
${\boldsymbol D}\cap\dom {\boldsymbol f}\neq\emp$. 
Appealing to \cite[Theorem~9.35 and Corollary~9.37]{Deut01}, we conclude that 
${\boldsymbol D}-\dom{\boldsymbol f}$ is a closed vector subspace.

\ref{p:sabang-mai2008iv}:
See \eqref{e:palawan-mai08-1}.

\ref{p:sabang-mai2008v}$\Rightarrow$\ref{p:sabang-mai2008iv}: 
See the proof of \cite[Theorem~6.3]{Baus93}.

\ref{p:sabang-mai2008vi}: Using
Lemma~\ref{l:palawan-mai08}\ref{l:palawan-mai08ii}\&\ref{l:palawan-mai08iii},
we obtain
${\boldsymbol 0}\in\sri{\boldsymbol C}$ $\Leftrightarrow$ 
${\boldsymbol 0}\in\sri({\boldsymbol D}-\dom {\boldsymbol f})$ 
$\Leftrightarrow$ 
${\boldsymbol 0}\in\reli({\boldsymbol D}-\dom {\boldsymbol f})$
$\Leftrightarrow$ 
${\boldsymbol 0}\in\reli{\boldsymbol D}-\reli\dom {\boldsymbol f}=
{\boldsymbol D}-\reli\dom {\boldsymbol f}$
$\Leftrightarrow$ $\boldsymbol D\cap\reli\dom {\boldsymbol f}\neq\emp$
$\Leftrightarrow$ $\bigcap_{i=1}^m\reli\dom f_i\neq\emp$.
\end{proof}

\section{Applications to signal and image processing}
\label{sec:4}

To illustrate the versatility of the proposed framework, we present 
three applications in signal and image processing. In each experiment, 
Algorithm~\ref{algo:1} is implemented with $\omega_i\equiv 1/m$, 
$\lambda_n\equiv 1.5$, and, since the proximity operators required by the 
algorithm will be computable in closed form, we can dispense with errors 
and set $a_{i,n}\equiv 0$ in \eqref{e:main1}. 
As a result, conditions~\ref{t:3iii} and~\ref{t:3iv} 
in Theorem~\ref{t:3} are straightforwardly satisfied. 
In each experiment, the number of iterations of the algorithm
is chosen large enough so that no significant improvement is gained 
by letting the algorithm run further.

\begin{figure}
\begin{center}
\includegraphics[width=9cm]{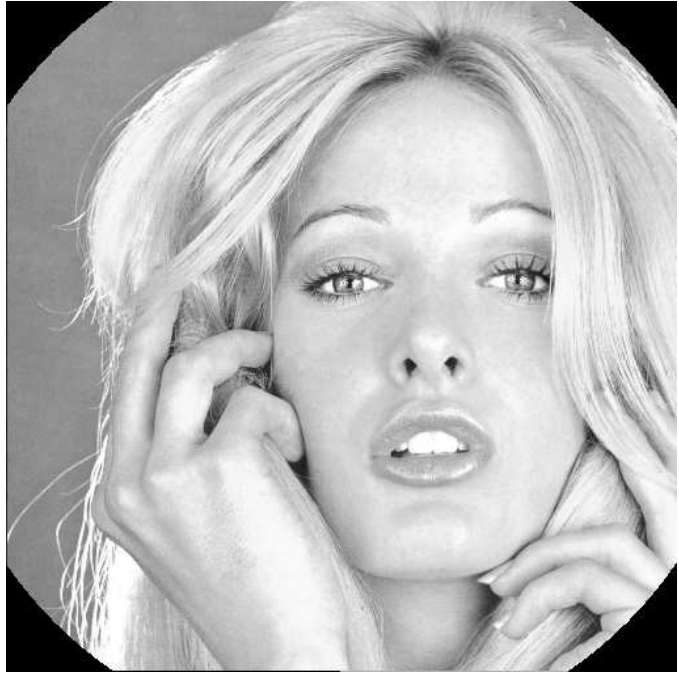}
\caption{Experiment 1. Original image.}
\label{fig:11}
\end{center}
\end{figure}

\begin{figure}
\begin{center}
\includegraphics[width=9cm]{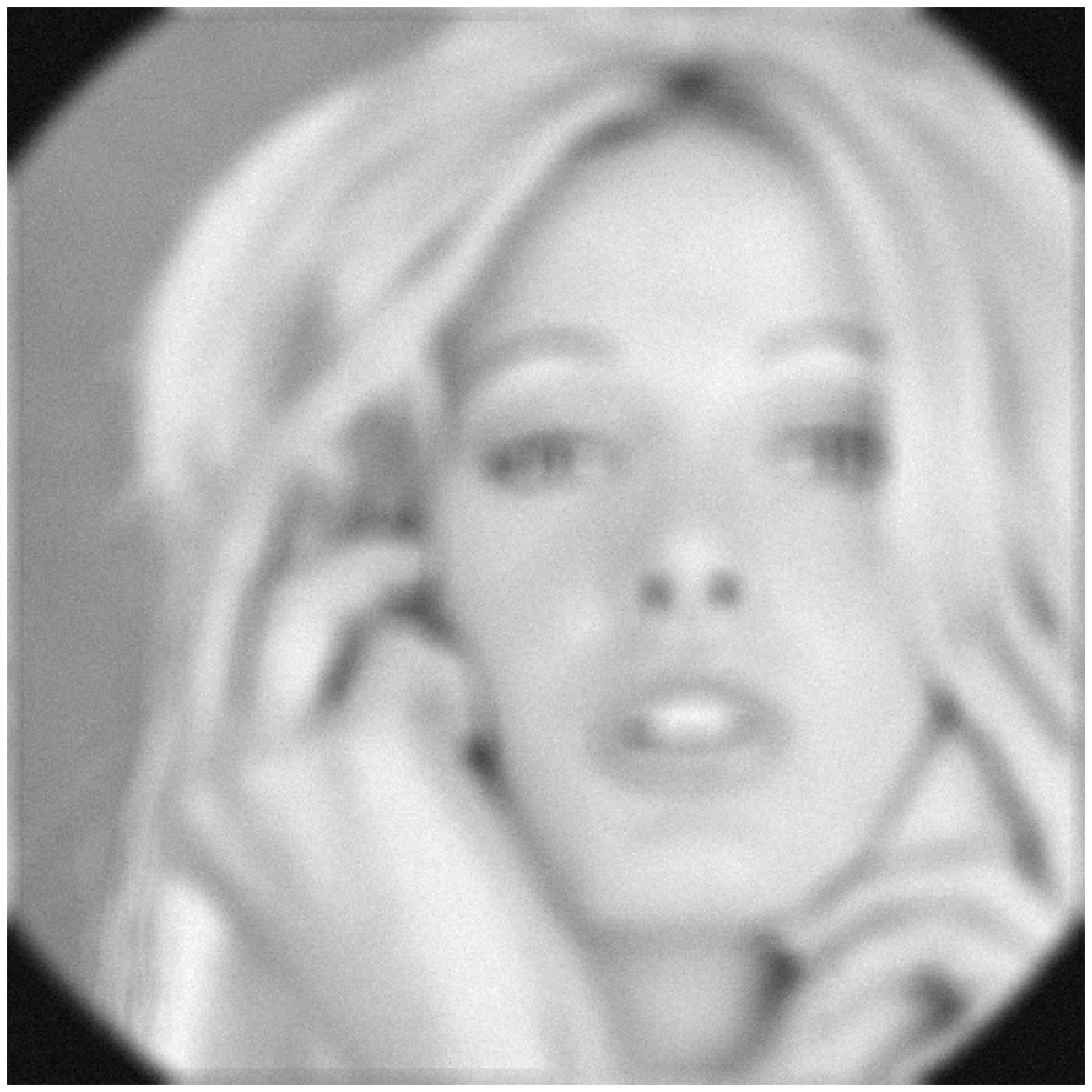}
\caption{Experiment 1. Degraded image.}
\label{fig:12}
\end{center}
\end{figure}

\begin{figure}
\begin{center}
\includegraphics[width=9cm]{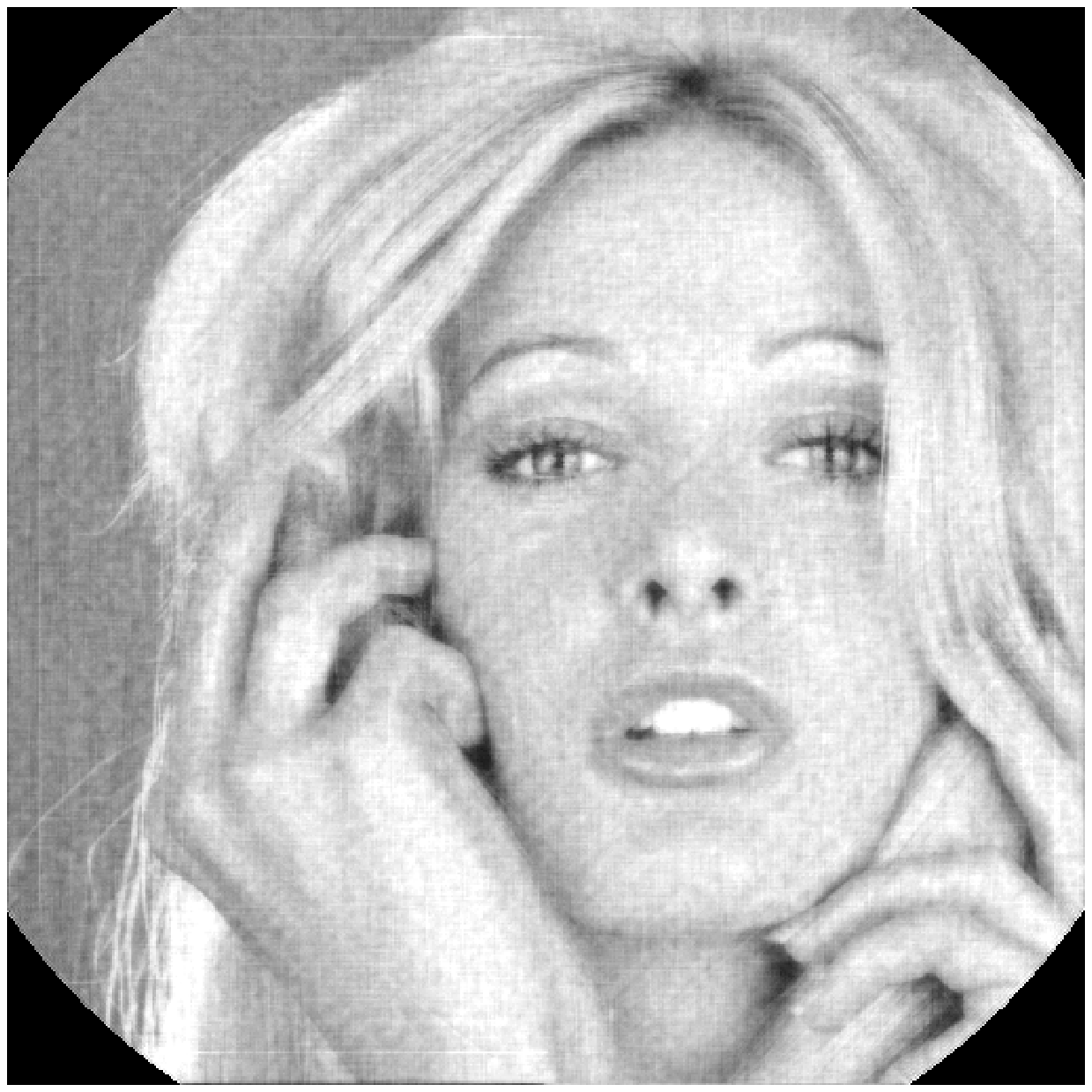}
\caption{Experiment 1. Image restored with 300 iterations of 
Algorithm \ref{algo:1}
($\gamma=1/4$).}
\label{fig:13}
\end{center}
\end{figure}

\subsection{Experiment 1}
\label{sec:exp1}
This first experiment is an image restoration problem in the standard 
Euclidean space $\HH=\RR^{N^2}$, where $N=512$. The original vignetted 
$N\times N$ image $\overline{x}$ is shown in figure~\ref{fig:11}
(the vignetting is modeled by a black area in the image corners). The
degraded image $z$ shown in figure~\ref{fig:12} is obtained via the
degradation model 
\begin{equation}
z=L\overline{x}+w,
\label{eq:addmodel}
\end{equation}
where $L$ is the two-dimensional convolution operator induced by a 
$15\times 15$ uniform kernel, and where $w$ is a realization of a 
zero-mean white Gaussian noise. The blurred image-to-noise ratio is 
$20\log_{10}(\|L\overline{x}\|/\|w\|)=31.75$~dB and the relative quadratic 
error with respect to the original image is 
$20\log_{10}(\|z-\overline{x}\|/\|\overline{x}\|)=-19.98$~dB.

The pixel values are known to fall in the interval $[0,255]$. In addition, 
the vignetting area $\mathbb S$ of the original image is known. This 
information leads to the constraint set
\begin{equation}
\label{e:14juillet2008-1}
C_1=[0,255]^{N^2}\cap\menge{x\in\HH}{x\,1_{\mathbb S}=\underline{0}},
\end{equation}
where $x1_{{\mathbb S}}$ denotes the coordinatewise multiplication 
of $x$ with the characteristic vector $1_{{\mathbb S}}$ of 
$\mathbb S$ (its $k$th coordinate is $1$ or $0$ according as
$k\in{\mathbb S}$ or $k\notin{\mathbb S}$), and where 
$\underline{0}$ the zero image. The mean value $\mu\in\left]0,255\right[$ 
of $\overline{x}$ is also known, which corresponds to the constraint set
\begin{equation}
\label{e:14juillet2008}
C_2=\menge{x\in\HH}{\scal{x}{\underline{1}}=N^2\mu},
\end{equation}
where $\underline{1}=[1,\ldots,1]^\top\in\RR^{N^2}$.
In addition, the phase of the discrete Fourier transform of the 
original image is measured over some frequency range 
$\mathbb{D}\subset\{0,\ldots,N^2-1\}$ \cite{Sign94,Levi83,Pora99}. 
If we denote by 
$\widehat{x}=\big(|\chi_k|\exp(i\angle\chi_k)\big)_{0\leq k\leq N^2-1}$
the discrete Fourier transform of an image $x\in\HH$ and by 
$(\phi_k)_{k\in\mathbb{D}}$ the known phase values, we obtain the 
constraint set 
\begin{equation}
C_3=\menge{x\in \HH}{(\forall k\in\mathbb{D})\;\angle\chi_k=\phi_k}.
\end{equation}
A constrained least-squares formulation of the problem is
\begin{equation}
\label{e:manille-16mai2008-11}
\underset{x\in C_1\cap C_2\cap C_3}{\mathrm{minimize}}\;\;\|Lx-z\|^2
\end{equation}
or, equivalently,
\begin{equation}
\label{e:manille-16mai2008-1}
\underset{x\in C_1\cap C_2}{\mathrm{minimize}}\;\;\iota_{C_3}(x)
+\|Lx-z\|^2.
\end{equation}
However, in most instances, the phase cannot be measured exactly. 
This is simulated by introducing a 5~\% perturbation 
on each of the phase components $(\phi_k)_{k\in\mathbb{D}}$. To take 
these uncertainties into account in \eqref{e:manille-16mai2008-1},
we replace the ``hard" potential $\iota_{C_3}$ by a smoothed version,
namely $\alpha d^p_{C_3}$, for some $\alpha\in\RPP$
and $p\in\left[1,\pinf\right[$.
This leads to the variational problem 
\begin{equation}
\label{e:manille-16mai2008-2}
\underset{x\in C_1\cap C_2}{\mathrm{minimize}}\;\;
\alpha d_{C_3}^p(x)+\|Lx-z\|^2,
\end{equation}
which is a special case of \eqref{e:prob1}, with
$m=4$, $f_1=\iota_{C_1}$, $f_2=\iota_{C_2}$,
$f_3=\alpha d^p_{C_3}$, and  $f_4=\|L\cdot-z\|^2$. 
Let us note that, since $C_1$ is bounded, 
condition~\ref{t:3i} in Theorem~\ref{t:3} is satisfied. In addition, it
follows from Proposition~\ref{p:sabang-mai2008}\ref{p:sabang-mai2008vi}
that condition~\ref{t:3ii} in Theorem~\ref{t:3} also holds. Indeed, set 
$E=\left]0,255\right[^{N^2}\cap A \cap C_2$, where 
$A=\menge{x\in\HH}{x1_{\mathbb S}=\underline{0}}$. Then it follows from
\eqref{e:14juillet2008} that
\begin{equation}
\frac{N^2\mu}{N^2-\card \mathbb{S}}\big(\underline{1}-1_{\mathbb S}\big)\in E. 
\end{equation}
Hence, since $A$ and $C_2$ are affine subspaces, \eqref{e:14juillet2008-1} 
and Lemma~\ref{l:palawan-mai08}\ref{l:palawan-mai08i} yield
\begin{equation}
\bigcap_{i=1}^4\reli\dom f_i
=\reli C_1\cap\reli C_2
=(\reli C_1)\cap C_2
=(\reli\left[0,255\right]^{N^2})\cap A\cap C_2
=E\neq\emp.
\end{equation}

Problem~\eqref{e:manille-16mai2008-2} is solved for the following scenario: 
$\mathbb D$ corresponds to a low frequency band including about 80~\% of 
the frequency components, $p=3/2$, and $\alpha = 10$. The proximity
operators required by Algorithm~\ref{algo:1} are obtained as follows.
First, $\prox_{f_1}$ and $\prox_{f_2}$ 
are respectively the projectors onto $C_1$ and $C_2$, which can be obtained 
explicitly \cite{Aiep96}. Next, $\prox_{f_3}$ is given in 
Example~\ref{ex:manila2008-05-07}. It involves $P_{C_3}$, which can be 
found in \cite{Aiep96}. Finally, $\prox_{f_4}$ is supplied by
Proposition~\ref{p:30avril2008}. Note that, since $L$ is a two-dimensional
convolutional blur, it can be approximated by a block circulant matrix and
hence \eqref{e:30avril2008} can be efficiently implemented in
the frequency domain via the fast Fourier transform \cite{Andr77}.
The restored image, shown in figure~\ref{fig:13}, is much sharper than 
the degraded image $z$ and it achieves a relative quadratic error 
of $-23.25$~dB with respect to the original image $\overline{x}$.

\subsection{Experiment 2}

In image recovery, variational formulations involving total variation 
\cite{Cham97,Rudi92,Stro03} or sparsity promoting potentials 
\cite{Bect04,Bred08,Cham98,Daub04} are popular. The objective of the 
present experiment is to show that it is possible to employ
more sophisticated, hybrid potentials.

\begin{figure}
\begin{center}
\includegraphics[width=8.5cm]{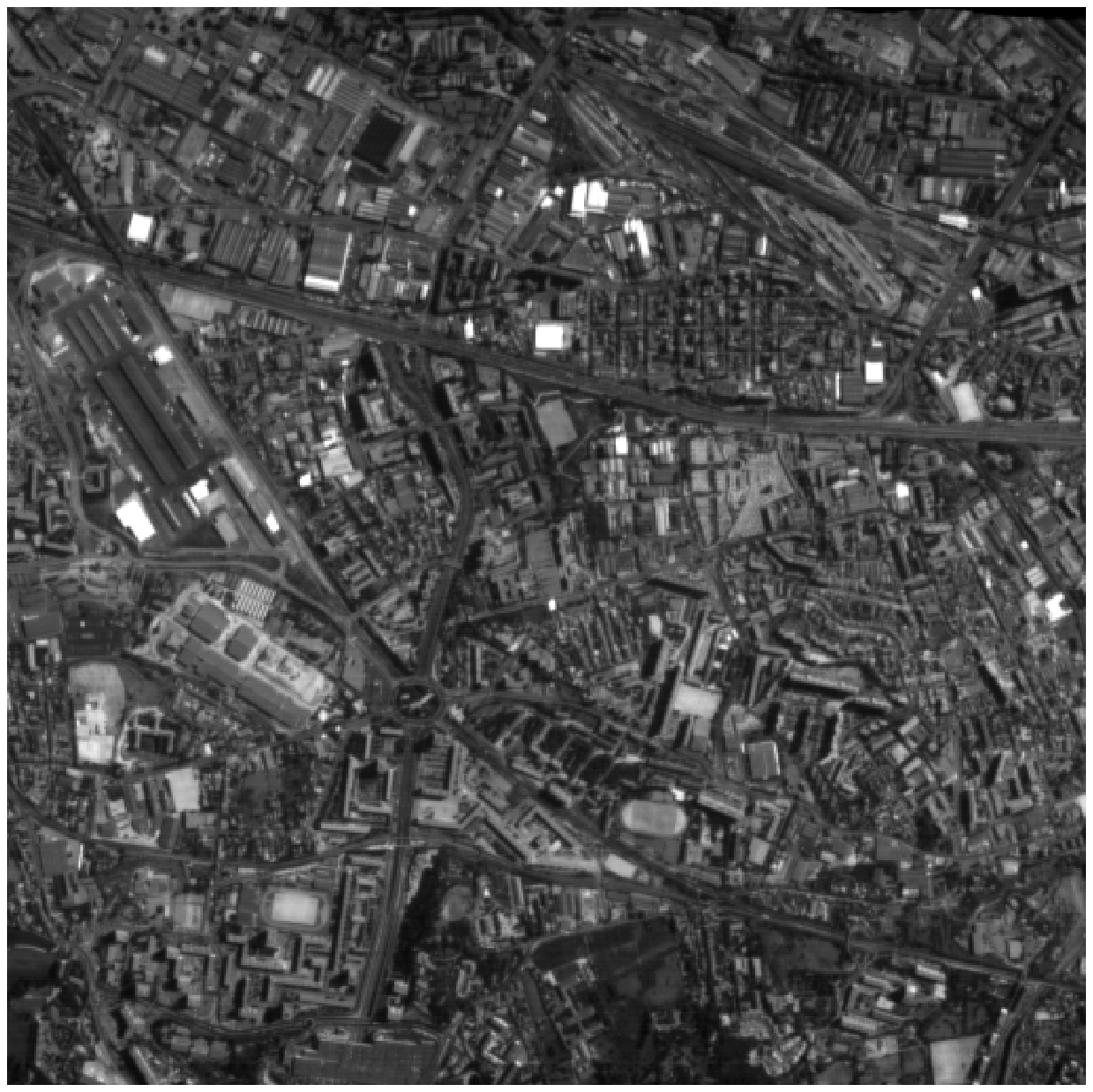}
\caption{Experiment 2. Original image.}
\label{fig:21}
\end{center}
\end{figure}

\begin{figure}
\begin{center}
\includegraphics[width=8.5cm]{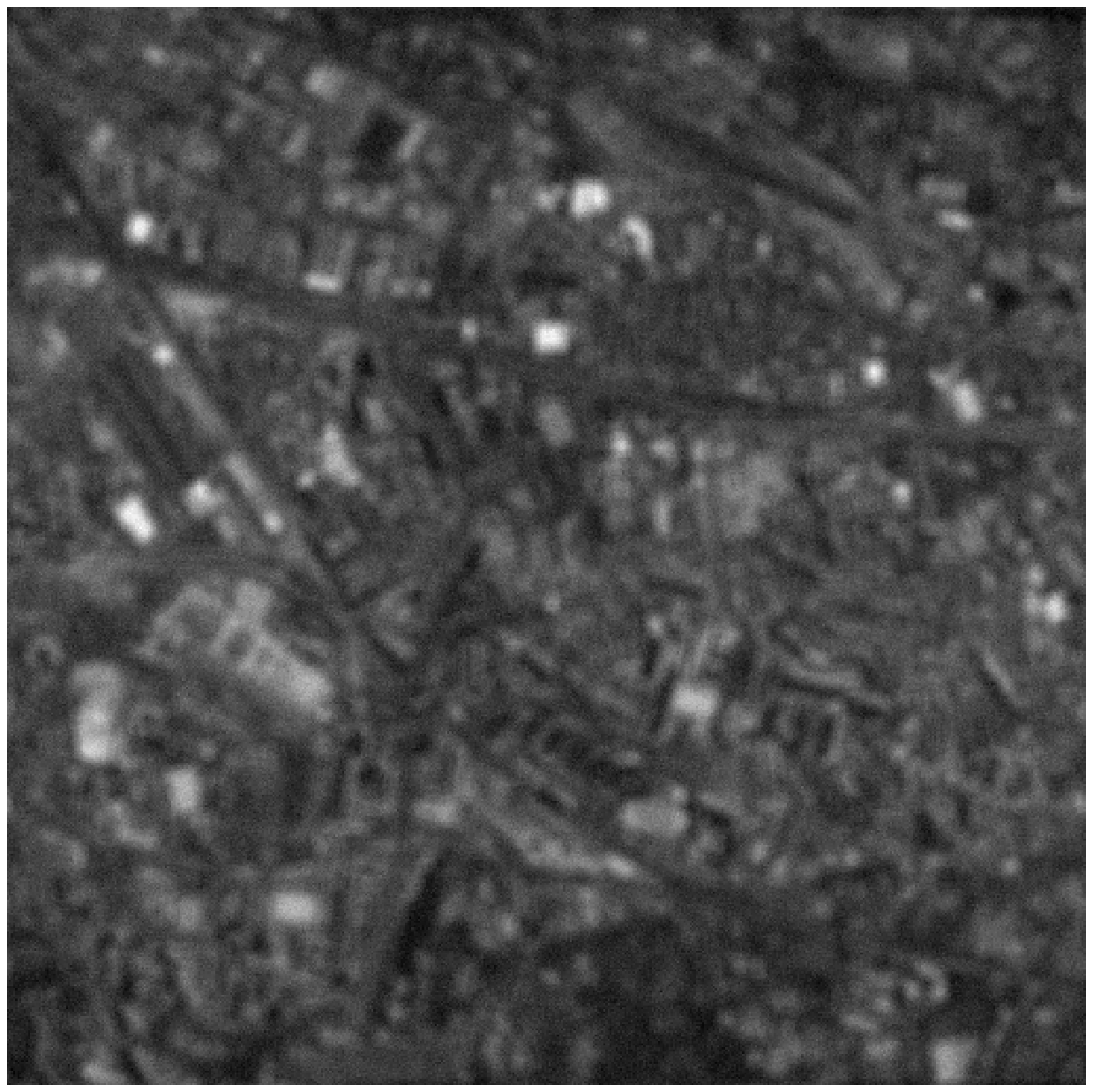}
\caption{Experiment 2. Degraded image.}
\label{fig:22}
\end{center}
\end{figure}

In order to simplify our presentation, we place ourselves in the Hilbert space 
$\GG$ of periodic discrete images $y=(\eta_{k,l})_{(k,l)\in\ZZ^2}$ with 
horizontal and vertical periods equal to $N$ ($N=512$), endowed with the 
standard Euclidean norm 
\begin{equation}
y\mapsto\sqrt{\sum_{k=0}^{N-1}\sum_{l=0}^{N-1}|\eta_{k,l}|^2}. 
\end{equation}
As usual, images of size $N\times N$ are viewed as elements of 
this space through periodization \cite{Andr77}.
The original 8-bit satellite image $\overline{y}\in \GG$ displayed in 
figure~\ref{fig:21} is degraded through the linear model
\begin{equation}
z=L\overline{y}+w,
\end{equation}
where $L$ is the two-dimensional periodic
convolution operator with a $7\times 7$ uniform kernel, and $w$ is a 
realization of a periodic zero-mean white Gaussian noise. 
The resulting degraded image $z\in\GG$ is shown in figure~\ref{fig:22}.
The blurred image-to-noise ratio is 
$20\log_{10}(\|L\overline{y}\|/\|w\|)=20.71~$dB 
and the relative quadratic error with respect to the original image is 
$20\log_{10}(\|z-\overline{y}\|/\|\overline{y}\|)=-12.02~$dB.

In the spirit of a number of recent investigations (see \cite{Chau07}
and the references therein), we use a tight frame representation of the 
images under consideration. This representation is defined through a 
synthesis operator $F^*$, which is a linear operator from $\HH=\RR^K$ 
to $\GG$ (with $K \geq N^2$)
such that
\begin{equation}
\label{e:tight}
F^*\circ F=\kappa \Id
\end{equation}
for some $\kappa\in\RPP$. Thus, the original image can be written as 
$\overline{y}=F^*\overline{x}$, where $\overline{x}\in \HH$ is a vector of
frame coefficients to be estimated. The rationale behind this approach is
that, by appropriately choosing the frame, a sparse representation
$\overline{x}$ of $\overline{y}$ can be achieved. 

The restoration problem is posed in the frame coefficient space $\HH$.
We use the constraint set imposing the range of the pixel values of 
the original image $\overline{y}$, namely
\begin{equation}
\label{e:C}
C=\menge{x\in\HH}{F^*x\in D},\;\:\text{where}\;\:
D=\menge{y\in\GG}{(\forall (k,l)\in\{0,\ldots,N-1\}^2)
\;\eta_{k,l}\in[0,255]},
\end{equation}
as well as three potentials. The first potential is the standard least-squares
data fidelity term $x\mapsto\|LF^*x-z\|^2$. The second potential
is the $\ell^1$ norm, which aims at promoting a sparse frame representation 
\cite{Chau07,Daub04,Trop06}. Finally, the third potential is the discrete 
total variation $\operatorname{tv}$, which aims at preserving piecewise
smooth areas and sharp edges \cite{Cham97,Rudi92,Stro03}. Using
the notation $(\eta_{k,l})_{(k,l)\in\ZZ^2}^\top=(\eta_{l,k})_{(k,l)\in\ZZ^2}$,
the discrete total variation of $y\in\GG$ is defined as
\begin{equation}
\label{e:8fyuihuio}
\mathrm{tv}(y)=\sum_{k=0}^{N-1}\sum_{l=0}^{N-1} 
\varrho_{k,l}\big(\nabla_{\!1}y,(\nabla_{\!1}(y^\top))^\top\big),
\end{equation}
where $\nabla_{\!1}\colon\GG\to\RR^{N\times N}$ is a discrete 
vertical gradient operator and where, for every 
$\{k,l,q,r\}\subset\{0,\ldots,N-1\}$, we set
\begin{equation}
\label{e:osaka2008-08-31}
\varrho_{k,l}=\varrho_{k,l,k,l}\:,
\end{equation}
with
\begin{equation}
\label{e:ricardo}
\varrho_{k,l,q,r}\colon\RR^{N\times N}\times\RR^{N\times N}\to\RR\colon
\Big(\big[\nu_{a,b}\big]_{0\leq a,b\leq N-1},
\big[\tilde{\nu}_{a,b}\big]_{0\leq a,b\leq N-1}\Big)
\mapsto\sqrt{|\nu_{k,l}|^2+|{\tilde{\nu}}_{q,r}|^2}.
\end{equation}
A common choice for the gradient operator is
$\nabla_{\!1}\colon y\mapsto[\eta_{k+1,l}-\eta_{k,l}]_{0\leq k,l\leq N-1}$.
As is customary in image processing \cite[Section~9.4]{Jain89}, 
we adopt here a horizontally smoothed version of this operator, namely,
\begin{align}
\label{eq:D1}
&\nabla_{\!1}\colon\GG\to\RR^{N\times N}\colon y
\mapsto\frac{1}{2}\big[\eta_{k+1,l+1}-\eta_{k,l+1}+
\eta_{k+1,l}-\eta_{k,l}\big]_{0\leq k,l\leq N-1}.
\end{align}
We thus arrive at a variational formulation of the form \eqref{e:prob1},
namely
\begin{equation}
\label{e:tubbs}
\underset{x\in\HH}{\mathrm{minimize}}\;\;\iota_C(x)+\|LF^*x-z\|^2
+\alpha\|x\|_{\ell^1}+\beta\mathrm{tv}(F^*x),
\end{equation}
where $\alpha$ and $\beta$ are in $\RPP$. Since $C$ is bounded, 
condition~\ref{t:3i} in Theorem~\ref{t:3} is satisfied. In addition, it
is clear from Proposition~\ref{p:sabang-mai2008}\ref{p:sabang-mai2008vi}
that condition~\ref{t:3ii} in Theorem~\ref{t:3} also holds. Indeed, all the
potentials in \eqref{e:tubbs} have full domain, except $\iota_C$. However,
Lemma~\ref{l:palawan-mai08}\ref{l:palawan-mai08ii} implies that
$\reli\dom\iota_C=\reli C\neq\emp$ since $\underline{0}\in C$.

\begin{figure}
\begin{center}
\includegraphics[width=8.5cm]{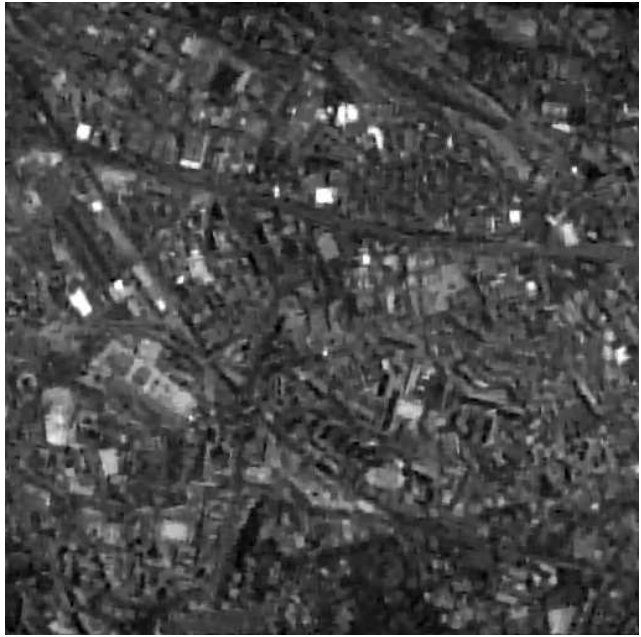}
\caption{Experiment 2. Image restored by \eqref{e:pbtv}, using 350 
iterations of Algorithm~\ref{algo:1} with $\gamma=150$.}
\label{fig:23}
\end{center}
\end{figure}

\begin{figure}
\begin{center}
\includegraphics[width=8.5cm]{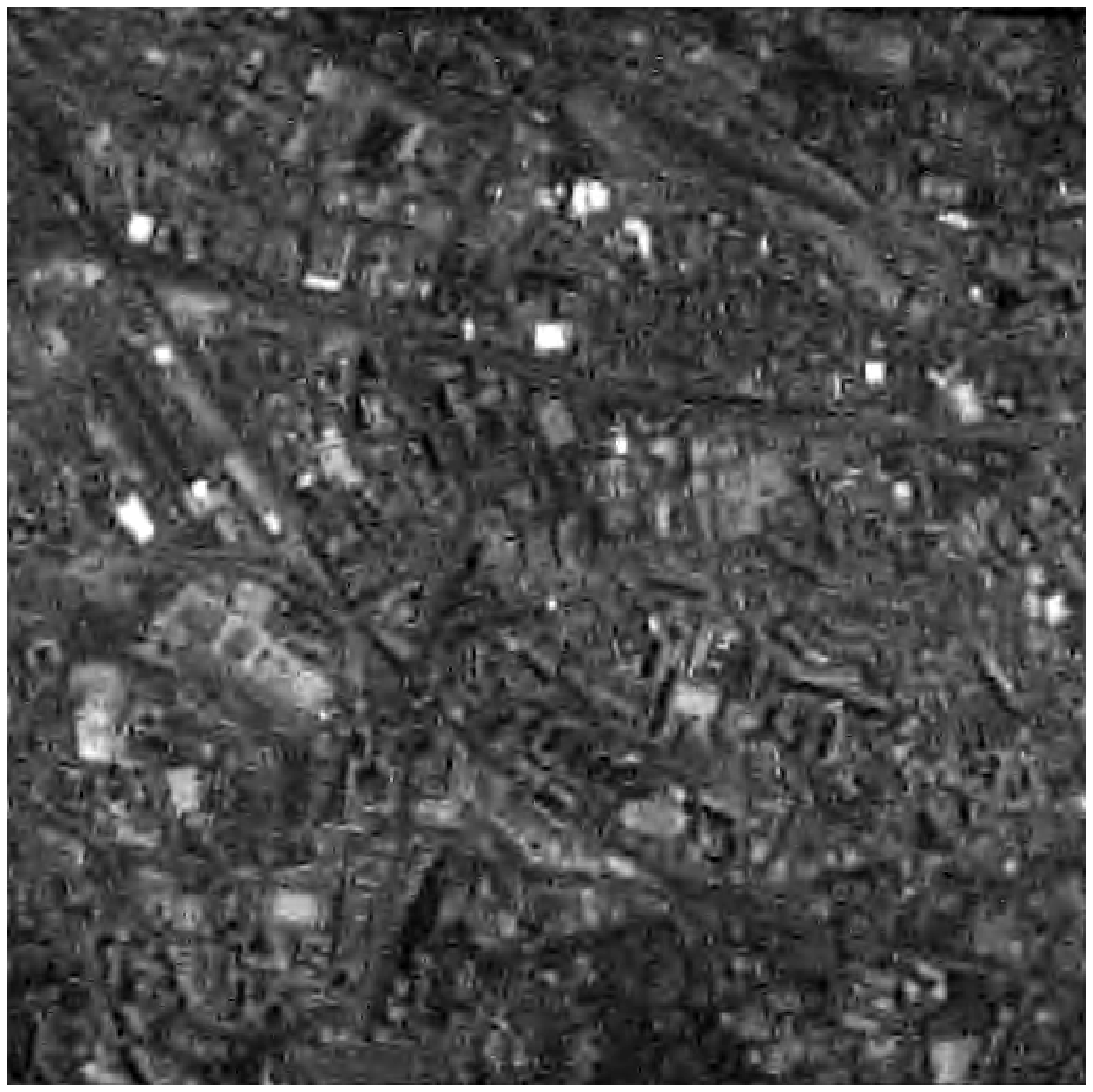}
\caption{Experiment 2. Image restored without the total variation potential
in \eqref{e:pbtv}, using 350 iterations of Algorithm~\ref{algo:1} 
with $\gamma=150$.}
\label{fig:24}
\end{center}
\end{figure}

\begin{figure}
\begin{center}
\includegraphics[width=8.5cm]{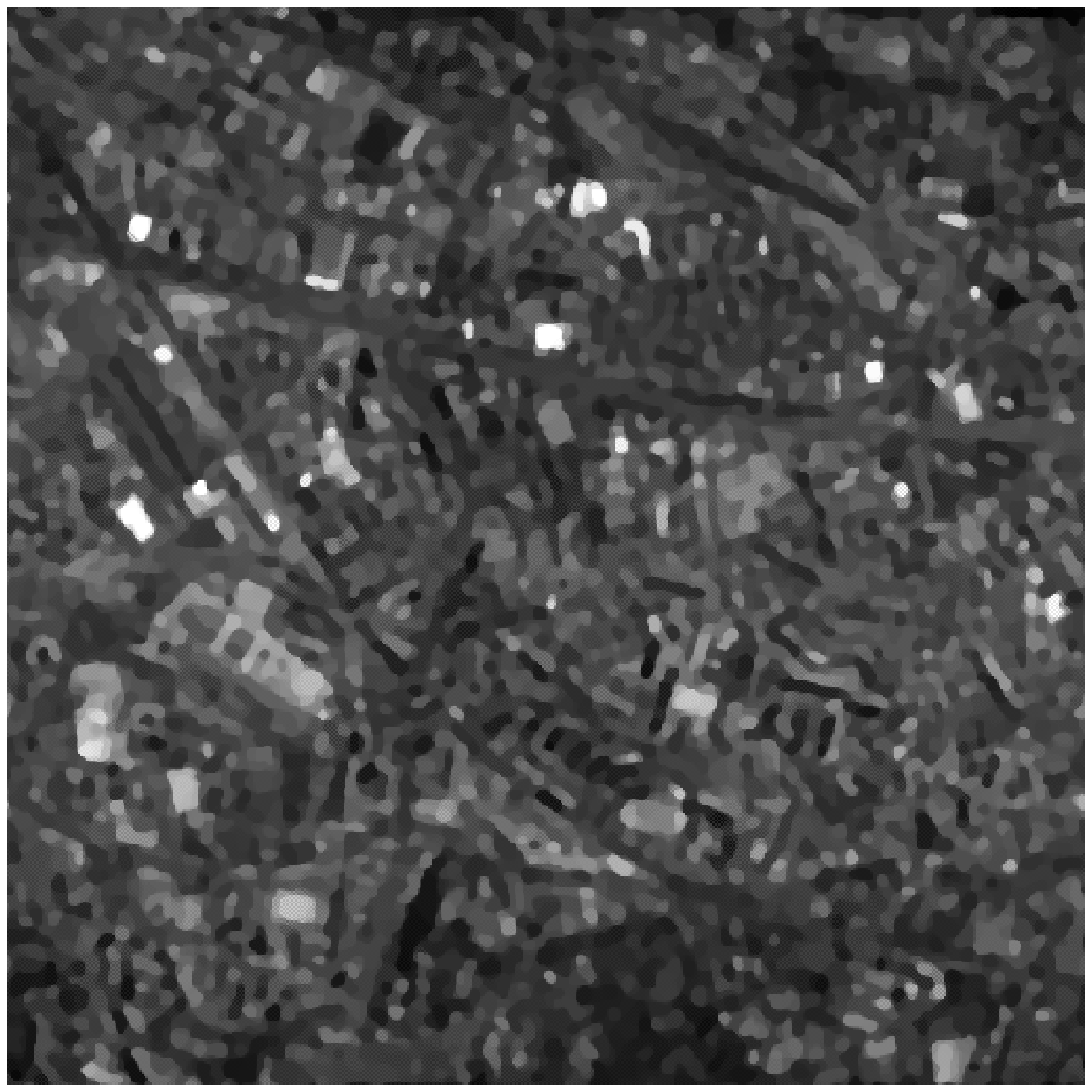}
\caption{Experiment 2. Image restored without the $\ell^1$ potential in
\eqref{e:pbtv}, using 350 iterations of Algorithm~\ref{algo:1} 
with $\gamma=150$.}
\label{fig:25}
\end{center}
\end{figure}

Although \eqref{e:tubbs} assumes the form of \eqref{e:prob1}, it is not 
directly exploitable by Algorithm~\ref{algo:1} because the proximity operator 
of $\operatorname{tv}\circ F^*$ cannot be computed explicitly. To circumvent 
this numerical hurdle, the total variation potential \eqref{e:8fyuihuio} is 
split in four terms and \eqref{e:tubbs} is rewritten as
\begin{equation}
\label{e:pbtv1}
\underset{x\in C}{\mathrm{minimize}}\;\;\|LF^*x-z\|^2
+\alpha\|x\|_{\ell^1}+\beta\sum_{i=0}^3\mathrm{tv}_i(F^*x),
\end{equation}
where
\begin{equation}
\label{e:15juillet2008}
(\forall (q,r) \in \{0,1\}^2)\quad
\mathrm{tv}_{q+2r}\colon\GG\to\RR\colon y\mapsto
\sum_{k=0}^{N/2-1}\sum_{l=0}^{N/2-1} 
\varrho_{2k+q,2l+r}\big(\nabla_{\!1}y,(\nabla_{\!1}(y^\top))^\top\big).
\end{equation}
For every $q$ and $r$ in $\{0,1\}$, let $\downarrow_{q,r}$ be the decimation 
operator given by
\begin{equation}
\downarrow_{q,r}\colon\RR^{2N\times 2N}\to\RR^{N\times N}\colon
v=\big[\nu_{k,l}\big]_{0\leq k,l\leq 2N-1}\mapsto 
\big[\nu_{2k+q,2l+r}\big]_{0\le k,l \leq N-1},
\end{equation}
and set 
\begin{align}
\label{e:15juillet2008-1}
U_{q+2r}\colon\GG\to\RR^{N\times N}\colon y
\mapsto\;
\downarrow_{q,r}
\begin{bmatrix}
\nabla_{\!0}y&
\nabla_{\!1}y\\
(\nabla_{\!1}(y^\top))^\top&
\nabla_{\!2}y
\end{bmatrix},
\end{align}
where $\nabla_1$ is defined in \eqref{eq:D1},
\begin{equation}
\label{eq:D0}
\nabla_0\colon\GG\to\RR^{N\times N}\colon y
\mapsto\frac{1}{2}\big[\eta_{k+1,l+1}+\eta_{k,l+1}
+\eta_{k+1,l}+\eta_{k,l}\big]_{0\leq k,l\leq N-1},
\end{equation}
and 
\begin{equation}
\label{eq:D3}
\nabla_{\!2}\colon\GG\to\RR^{N\times N}\colon y
\mapsto\frac{1}{2}\big[\eta_{k+1,l+1}-\eta_{k,l+1}
-\eta_{k+1,l} +\eta_{k,l}\big]_{0\leq k,l\leq N-1}.
\end{equation}
Moreover, set 
\begin{equation}
\label{e:16juillet2008}
h\colon\RR^{N\times N}\to\RR\colon v\mapsto
\sum_{k=0}^{N/2-1}\sum_{l=0}^{N/2-1}
\varrho_{k,l+N/2,k+N/2,l}(v,v).
\end{equation}
Then it follows from \eqref{e:15juillet2008} and 
\eqref{e:15juillet2008-1} that 
\begin{equation}
\label{e:15juillet2008-2}
(\forall i\in\{0,1,2,3\})\quad\mathrm{tv}_i=h\circ U_i.
\end{equation}
Hence, \eqref{e:pbtv1} becomes
\begin{equation}
\label{e:pbtv}
\underset{x\in C}{\mathrm{minimize}}\;\;\|LF^*x-z\|^2
+\alpha \|x\|_{\ell^1}+\beta\sum_{i=0}^3h(U_iF^*x).
\end{equation}
Problem~\eqref{e:pbtv} is a specialization of \eqref{e:prob1}, in 
which $m=7$, $f_1=\iota_C$, $f_2=\|LF^*\cdot-z\|^2$, 
$f_3=\alpha\|\cdot\|_{\ell^1}$, and 
$f_{i+4}=\beta\,h\circ U_i\circ F^*$ for $i\in\{0,1,2,3\}$.
To implement Algorithm~\ref{algo:1}, we need the expressions of
the proximity operators of these functions. The proximity operator of 
$f_1$ can be calculated by first observing that the projection onto 
the set $D$ of \eqref{e:C} is explicit, and by then applying 
Lemma~\ref{l:1938}, which states that \eqref{e:tight} and \eqref{e:C} 
imply that
\begin{equation}
\prox_{f_1}=\prox_{\iota_D\circ F^*}
=\Id+\frac{1}{\kappa}F\circ(\prox_{\iota_D}-\Id)\circ F^*
=\Id+\frac{1}{\kappa}F\circ(P_D-\Id)\circ F^*.
\end{equation}
On the other hand, the proximity operator of $f_2$ can be derived from
Proposition~\ref{p:30avril2008} using a frequency domain implementation 
(as in section~\ref{sec:exp1}), 
and by again invoking Lemma~\ref{l:1938}. Next, the proximity operator 
of $f_3$ can be found in \cite[Example~2.20]{Smms05}. Finally, the operators 
$(\prox_{f_i})_{4\leq i\leq 7}$ are provided by the following fact.

\begin{proposition}
Set $\Pi\colon\RR^{N\times N}\to\RR^{N\times N}\colon 
v=\big[\nu_{k,l}\big]_{0\leq k,l\leq N-1}\mapsto
\big[\pi_{k,l}\big]_{0\leq k,l\leq N-1}$, where
\begin{multline}
\label{e:fine}
\big(\forall (k,l)\in\{0,\ldots,N/2-1\}^2\big)\quad
\begin{cases}
\pi_{k,l}=\nu_{k,l}\\
\pi_{k+N/2,l+N/2}=\nu_{k+N/2,l+N/2}\\
\pi_{k,l+N/2}=\sigma_{k,l}(v)\,\nu_{k,l+N/2}\\
\pi_{k+N/2,l}=\sigma_{k,l}(v)\,\nu_{k+N/2,l}
\end{cases}
\quad\text{with}\\[3mm]
\sigma_{k,l}\colon v\mapsto
\begin{cases}
1-\displaystyle{\frac{\kappa\beta}{\sqrt{|\nu_{k,l+N/2}|^2+|\nu_{k+N/2,l}|^2}}},
&\text{if}\;\;\sqrt{|\nu_{k,l+N/2}|^2+|\nu_{k+N/2,l}|^2}\geq\kappa\beta;\\
0,&\text{otherwise.}
\end{cases}
\end{multline}
Then, for every $i\in\{0,1,2,3\}$,
\begin{equation}
\label{eq:proxtvj}
\prox_{f_{i+4}}=\Id+\frac{1}{\kappa} F
\circ (U_i^*\circ\Pi\circ U_i-\Id)\circ F^*.
\end{equation}
\end{proposition}
\begin{proof}
Set $\varphi\colon\RR^2\to\RR\colon(\xi_1,\xi_2)
\mapsto\kappa\beta\sqrt{|\xi_1|^2+|\xi_2|^2}$. By applying 
Proposition~\ref{p:manille2008-05-05}\ref{p:manille2008-05-05i} in 
$\RR^2$ with the set $\{(0,0)\}$, we obtain
\begin{equation}
\label{e:castillo}
(\forall (\xi_1,\xi_2)\in\RR^2)\quad
\prox_{\varphi}(\xi_1,\xi_2)=
\begin{cases}
\bigg(1-\Frac{\kappa\beta}{\sqrt{|\xi_1|^2+|\xi_2|^2}}\bigg)(\xi_1,\xi_2),
&\text{if}\;\sqrt{|\xi_1|^2+|\xi_2|^2}\geq\kappa\beta;\\
0,&\text{otherwise.}
\end{cases}
\end{equation}
Now set $p=[\pi_{k,l}]_{0\leq k,l\leq N-1}=\prox_{\kappa\beta h}v$.
In view of \eqref{e:prox}, \eqref{e:16juillet2008}, and
\eqref{e:ricardo}, $p$ minimizes over $\tilde p\in\RR^{N\times N}$ 
the cost
\begin{align}
\kappa\beta h(\tilde p)+\frac{1}{2}\|v-\tilde p\|^2
&=\kappa\beta\sum_{k=0}^{N/2-1}\sum_{l=0}^{N/2-1}
\varrho_{k,l+N/2,k+N/2,l}(\tilde p,\tilde p)+\frac{1}{2}
\sum_{k=0}^{N-1}\sum_{l=0}^{N-1}|\nu_{k,l}-\tilde \pi_{k,l}|^2\nonumber\\
&=\sum_{k=0}^{N/2-1}\sum_{l=0}^{N/2-1}\Big(\kappa\beta
\sqrt{|\tilde\pi_{k,l+N/2}|^2+|\tilde\pi_{k+N/2,l}|^2}\nonumber\\
&\hskip 34mm\;+\frac{1}{2}\big(|\nu_{k,l+N/2}-\tilde\pi_{k,l+N/2}|^2
+|\nu_{k+N/2,l}-\tilde\pi_{k+N/2,l}|^2\big)\Big)\nonumber\\
&\quad\;+\frac{1}{2}\sum_{k=0}^{N/2-1}\sum_{l=0}^{N/2-1}\big(|\nu_{k,l}-
\tilde\pi_{k,l}|^2+|\nu_{k+N/2,l+N/2}-\tilde\pi_{k+N/2,l+N/2}|^2\big).
\end{align}
Therefore,
\begin{equation}
(\forall (k,l)\in\{0,\ldots,N/2-1\}^2)\quad
\begin{cases}
(\pi_{k,l+N/2},\pi_{k+N/2,l})=\prox_\varphi(\nu_{k,l+N/2},\nu_{k+N/2,l}),\\
\pi_{k,l}=\nu_{k,l},\\
\pi_{k+N/2,l+N/2}=\nu_{k+N/2,l+N/2}.
\end{cases}
\end{equation}
Appealing to \eqref{e:fine} and \eqref{e:castillo}, we obtain 
$\Pi=\prox_{\kappa\beta h}$. Now, let $i\in\{0,1,2,3\}$. It follows from 
\eqref{e:15juillet2008-1} that $U_i$ is a separable two-dimensional 
Haar-like orthogonal operator \cite[Section~5.9]{Jain89}. Hence, 
appealing to \eqref{e:tight}, we obtain 
$(U_i\circ F^*)\circ (U_i \circ F^*)^*=\kappa\Id$.
In turn, Lemma~\ref{l:1938} yields 
\begin{align}
\prox_{f_{i+4}}&=\prox_{\beta\,h\circ(U_i\circ F^*)}\nonumber\\
&=\Id+\frac{1}{\kappa} 
(U_i\circ F^*)^*\circ(\prox_{\kappa\beta h}-\Id)\circ(U_i\circ F^*)\nonumber\\
&=\Id+\frac{1}{\kappa} 
F\circ (U_i^*\circ\Pi\circ\, U_i-\Id)\circ F^*,
\end{align} 
which completes the proof.
\end{proof}

In \eqref{e:pbtv}, we employ a tight frame ($\kappa=4$) resulting from the 
concatenation of four shifted separable dyadic orthonormal wavelet 
decompositions \cite{Pesq96} carried out over 4 resolution levels. 
The shift parameters are $(0,0)$, $(1,0)$, $(0,1)$, and $(1,1)$. 
In addition, symlet filters \cite{Daub92} of length 8 are used.
The parameters $\alpha$ and $\beta$ have been adjusted so as to minimize
the error with respect to the original image $\overline{y}$.
%alpha = 0.008 beta = 0.03
The restored image we obtain is displayed in 
figure~\ref{fig:23}. It achieves a relative mean-square error with respect 
to $\overline{y}$ of $-14.82$~dB. 
For comparison, the result obtained without the total variation potential
in \eqref{e:pbtv} is shown in figure~\ref{fig:24} (error of $-14.06$~dB), 
%alpha = 0.03 beta = 0
and the result obtained without the $\ell^1$ potential in \eqref{e:pbtv} 
is shown in figure~\ref{fig:25} (error of $-13.70$~dB).
%alpha = 0 beta = 0.1
It can be observed that the image of figure~\ref{fig:24} 
suffers from small visual artifacts, whereas the details in 
figure~\ref{fig:25} are not sharp. This shows the 
advantage of combining an $\ell^1$ potential and a total variation
potential.

\subsection{Experiment 3}

We revisit via the variational formulation \eqref{e:prob1} a pulse 
shape design problem investigated in \cite{Sign99} in a more 
restrictive setting (see also \cite{Noba95} for the original two-constraint
formulation). This problem illustrates further ramifications
of the proposed algorithm. 

The problem is to design a pulse shape for digital communications. The 
signal space is the standard Euclidean space $\HH=\RR^N$, where $N=1024$ 
is the number of samples of the discrete pulse (the underlying 
sampling rate is 2560~Hz). Five constraints arise from engineering 
specifications. We denote by $x=(\xi_k)_{0 \leq k \leq N-1}$ a signal 
in $\HH$ and by $\widehat{x}=(\chi_k)_{0 \leq k\leq N-1}$ its 
discrete Fourier transform.
\begin{itemize}
\item 
The Fourier transform of the pulse should vanish at the zero frequency and
at integer multiples of 50~Hz. This constraint is associated with the set
\begin{equation}
C_1=\menge{x\in\HH}{\widehat{x}\,1_{\mathbb{D}_1}=\underline{0}},
\end{equation}
where $\mathbb{D}_1$ is the set of discrete frequencies at which 
$\widehat{x}$ should vanish.
\item 
The modulus of the Fourier transform of the pulse should no exceed 
a prescribed bound $\rho>0$ beyond 300~Hz. This constraint is associated 
with the set
\begin{equation}
C_2=\menge{x\in\HH}{(\forall k\in\mathbb{D}_2)\;|\chi_k|\leq\rho},
\end{equation}
where $\mathbb{D}_2$ represents frequencies beyond 300~Hz.
\item 
The energy of the pulse should not exceed a prescribed bound $\mu^2>0$ in 
order not to interfere with other systems. The associated set is
\begin{equation}
C_3=\menge{x\in\HH}{\|x\|\leq\mu}.
\end{equation}
\item 
The pulse should be symmetric about its mid-point, where its value should
be equal to $1$. This corresponds to the set
\begin{equation}
C_4=\menge{x\in\HH}{\xi_{N/2} =1\;\text{and}\;
(\forall k \in \{0,\ldots,N/2\})\;\xi_k=\xi_{N-1-k}}.
\end{equation}
\item 
The duration of the pulse should be 50~ms and it should have periodic zero
crossings every 3.125~ms. This leads to the set
\begin{equation}
C_5=\menge{x\in\HH}{x\,1_{\mathbb S}=\underline{0}},
\end{equation}
where $\mathbb{S}$ is the set of time indices in the zero areas.
\end{itemize}

\begin{figure}
\begin{center}
\includegraphics[width=9cm]{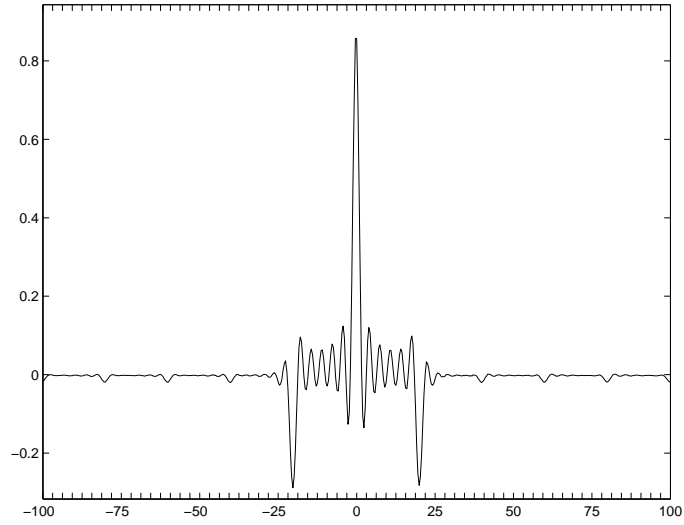}
\caption{Experiment 3. Pulse (amplitude versus time in ms) synthesized 
using 100 iterations of Algorithm~\ref{algo:1} with $\gamma=1/5$.}
\label{fig:31}
\end{center}
\end{figure}

\begin{figure}
\begin{center}
\includegraphics[width=9cm]{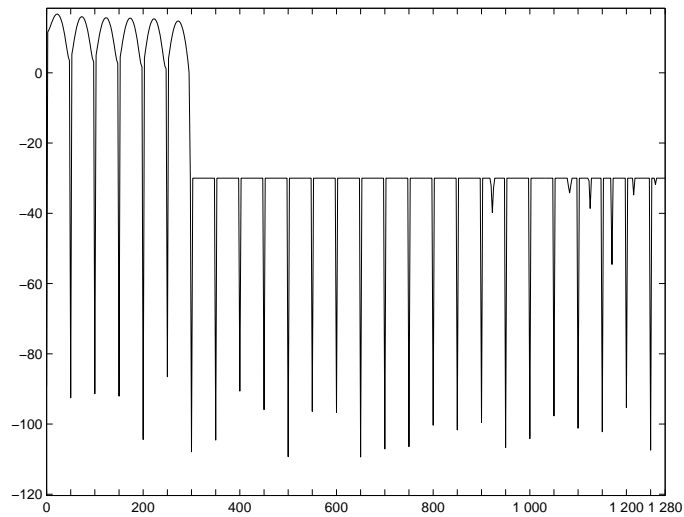}
\caption{Experiment 3. Fourier transform (amplitude in dB versus 
frequency in Hz) of the pulse of figure~\ref{fig:31}.} 
\label{fig:32}
\end{center}
\end{figure}

In this problem, $C_1$, $C_2$, and $C_3$ are hard constraints that must be
satisfied, whereas the other constraints are soft ones that are
incorporated via powers of distance potentials. This leads to the
variational formulation
\begin{equation}
\label{e:scenar1}
\underset{x\in C_1\cap C_2 \cap C_3}{\mathrm{minimize}}\;\;
d_{C_4}^{p_4}(x)+d_{C_5}^{p_5}(x),
\end{equation}
where $p_4$ and $p_5$ are in $[1,\pinf[$.
The design problem is thus cast in the general form of
\eqref{e:prob1}, with $m=5$,
$f_i=\iota_{C_i}$ for $i\in\{1,2,3\}$, 
and $f_i=d_{C_i}^{p_i}$ for $i \in \{4,5\}$.
Since $C_3$ is bounded, 
condition~\ref{t:3i} in Theorem~\ref{t:3} holds. In addition, it follows
from Proposition~\ref{p:sabang-mai2008}\ref{p:sabang-mai2008vi} that
condition~\ref{t:3ii} in Theorem~\ref{t:3} is satisfied. 
Indeed, 
\begin{equation}
0\in C_1\cap\menge{x\in \HH}
{(\forall k\in\mathbb{D}_2)\;|\chi_k|<\rho}\cap\menge{x\in\HH}{\|x\|<\mu}=
\bigcap_{i=1}^5\reli\dom f_i.
\end{equation}
Let us emphasize that our approach is applicable to any value of 
$(p_4,p_5)\in[1,\pinf[^2$. The proximity operators of $f_4$ and 
$f_5$ are supplied by Proposition~\ref{p:manille2008-05-05}, 
whereas the other proximity operators are the projectors onto 
$(C_i)_{1\leq i\leq 3}$, which are straightforward \cite{Sign99}.
A solution to \eqref{e:scenar1} when $p_4=p_5=2$, $\rho=10^{-3/2}$, and 
$\mu=2$ is shown in figure~\ref{fig:31} and its Fourier transform is 
shown in figure~\ref{fig:32}.
As is apparent in figure~\ref{fig:31}, the constraints corresponding to
$C_4$ and $C_5$ are not satisfied. Forcing $C_4\cap C_5$ as a hard
constraint would therefore result in an infeasible problem. Finally, 
figure~\ref{fig:32} shows that $C_2$ induces a 30~dB attenuation in the 
stop-band (beyond 300~Hz), in agreement with the value chosen for $\rho$.

\section{Concluding remarks}
\label{sec:5}

We have proposed a proximal method for solving inverse problems that 
can be decomposed into the minimization of a sum of lower
semicontinuous convex potentials. The algorithms currently in 
use in inverse problems are restricted to at most two nonsmooth 
potentials, which excludes many important scenarios and offers limited
flexibility in terms of numerical implementation. By contrast, 
the algorithm proposed in the paper can handle an arbitrary number of 
nonsmooth potentials. It involves each potential by means 
of its own proximity operator, and activates these
operators in parallel at each iteration. The versatility
of the method is demonstrated through applications in signal and image
recovery that illustrate various decomposition schemes, including one in
which total variation is mixed up with other nonsmooth potentials.
%In our future work, we will investigate conditions for the strong convergence
%of the proposed algorithm.

\end{document}